\documentclass[prd,floatfix]{revtex4}

\newcommand{\BA}{\begin{array}}
\newcommand{\EA}{\end{array}}
\newcommand{\BE}{\begin{equation}}
\newcommand{\EE}{\end{equation}}
\newcommand{\BEA}{\begin{eqnarray}}
\newcommand{\EEA}{\end{eqnarray}}
\newcommand{\BM}{\begin{bmatrix}}
\newcommand{\EM}{\end{bmatrix}}
\newcommand{\bdm}{\begin{displaymath}}
\newcommand{\edm}{\end{displaymath}}

\usepackage{amssymb, amsfonts, amsmath, graphicx}
\usepackage{grffile}
\usepackage{dcolumn}
\usepackage{bm}
\usepackage{epsfig}
\usepackage{mathrsfs}  
\usepackage{subfigure}
\usepackage{multirow}
\usepackage{epstopdf}

\begin{document}

\title{A Study of Space-Time Discretizations for the Dirac Equation}

\author{Robert Vaselaar$^{1}$, Hyun Lim$^{1}$, Jung-Han Kimn$^{1}$}

\affiliation{$1$ Department of Mathematics and Statistics, South Dakota State University, Brookings, SD 57007}

\begin{abstract}
We study several numerical discretization techniques for the one-space plus one-time dimensional Dirac equation, including finite difference and space-time finite element methods. Two finite difference schemes and several space-time finite elements function spaces are analyzed with respect to known analytic solutions. Further we propose a finite element discretization along the equations' characteristic lines, creating diamond-shaped elements in the space-time plane. We show that the diamond shaped elements allow for physically intuitive boundary conditions, improve numerical efficiency, and reduce the overall error of the computed solution as compared to the other finite difference and space-time finite element discretizations studied in this paper.
\end{abstract}
\maketitle

\section{Introduction}
\label{sec:dirac}

The Dirac equation governs all spin$-\frac 1 2$ particles, known as fermions.
While solutions to the Dirac equation  may be used to derive quantifiable predictions of particle behavior from quantum physics, these solutions are sometimes difficult to find in experimentally interesting scenarios. 
Numerical methods for the Dirac equation may be able to bridge some of the gap between theoretical and experimental particle physics \cite{fillion2012, widom1996neutrino}. 

These include many finite difference based methods for lattice quantum chromodynamics, such as the Dirac Wilson equation \cite{PhysRevD.10.2445} which have been used in conjunction with modern numerical methods such as Krylov subspace solvers \cite{Sakurai2010113, Nakamura201234} and scalable additive Schwarz preconditioners \cite{Luscher2004209}.
The limitations of these methods are also well known and include the inability to reconcile all limitations simultaneously. This is particularly important when considering the problem of fermion doubling, a condition where the number of particles considered must naturally double for each space-time dimension included on the lattice, and chiral symmetry, which is usually broken by most numerical methods that prevent fermion doubling \cite{Nielsen1981219, Chandrasekharan2004373}.

Other numerical methods for the Dirac equation include radial formulations created to investigate the energy spectrum of heavy atomic ions \cite{almanasrehgconverg, kullie2004, kullie2001, Desclaux2003453}. These are based on the Dirac equation in the presence of a coulomb potential resulting a relativistic eigenvalue problem, using both finite difference and finite element numerical methods.

The finite element method has also been used to calculate the propagation of free fermions in space. Analysis of the finite element method combined with Crank-Nicholson time stepping scheme demonstrates that solutions may show inconsistent and impossible physical behavior, such as superluminal propagation, depending on the step size and propagation method used \cite{muller1998}. 
Since using Lagrangian interpolation elements in one dimension are algebraically similar to using finite differences, it is also natural that this choice of function space has the same problems of fermion doubling and numerical instability as its finite difference relative\cite{muller1998}. 
In this implementation, physically consistent behavior of the particle depended on particle momentum, finite element size, and time step size chosen. 


In this paper, several implicit space-time discretizations based on the finite difference and Galerkin methods are presented. This presentation will show that simulation behavior is directly affected by choice of discretization method and function space. 
The problem domain is then rotated by $45^\circ$ in the space-time plane, forming diamond-shaped tensor elements, and the solution is recalculated using the rotated domain. 
This rotated domain shows substantially reduced error and improved performance when compared to the other space-time discretizations listed here.
The goal of this research is to create a discrete form of the Dirac equation that shows good agreement with the analytic solution as well as low error and the absence of faster-than-light propagation. 
Further, we would prefer a solution that does not modify the original Dirac operator in order keep as many of its original physical properties as possible.


This paper is organized as following;
First we present the weak form of the gauge-free Dirac Equation in section \ref{sec:STFEM}.  In section  \ref{sec:Stad_Numer} two space-time finite difference and one finite element method are presented along with their numerical results to observe their behavior and performance. Then in section \ref{sec:TEBA} three finite element discretizations using space-time tensor elements are presented along with their numerical results. Sources of possible simulation error are also presented and analyzed.  Proceeding from the discussion of error we propose our diamond-shaped approach in section \ref{sec:Diamond} and show how this approach addresses the errors observed and improves simulation efficiency. 
We conclude by discussing future research opportunities in Section 
\ref{sec:Conc}.

\section{  Space-Time Methods}
\label{sec:STFEM}

\subsection{Weak Formulation of the Dirac Equation}

The one dimensional Dirac operator may be expressed as follows
\begin{eqnarray}
\widehat D = \left(-i\hbar I\partial_t - i\hbar c \sigma_1\partial_x + mc^2\sigma_0 \right) \label{eqn:dirac1dOp}
\end{eqnarray}

where $\sigma_0, \sigma_1$ are the usual Pauli matrices defined as

\begin{equation}\label{eqn:sig0}
 \sigma_0 = \begin{bmatrix}
1 & 0\\ 
0 &  -1
\end{bmatrix}
\end{equation}

\begin{equation}\label{eqn:sig1}
 \sigma_1 = \begin{bmatrix}
0 & 1\\ 
1 &  0
\end{bmatrix} 
\end{equation}

Here the Pauli matrices are chosen such that the variables $x$ and $t$ form a Minkowski space-time, which an essential relationship in the Dirac equation.

In this case we will consider the initial value problem given by

\begin{eqnarray}
\widehat D  \widehat \Psi = 0 \label{eqn:dirac1d} \hspace{0.15cm} {\rm on} \hspace{0.15cm} \Omega \times [0,T]\\
\widehat \Psi (\cdot ,0) = \widehat \Psi^0 \label{eqn:psi0}
\end{eqnarray}

In the gage-free case analytic solutions may be computed directly which will give us a basis for comparison
for our numerical results. 

Using the continuous time Galerkin method the weak form may be expressed as follows. 
The objective is to find $\Phi,\Psi \in H^1_0 \left( \Omega \times [0,T] \right ) $ such that
\begin{equation}
\label{eqn:weakform}
\int_{\Omega \times [0,T]}  \widehat \Psi^*\widehat D \widehat \Psi ds = 0
\end{equation}

\section{Numerical Results of Different  Numerical Approaches}
\label{sec:Stad_Numer}

\begin{figure}[h]
\centering
\subfigure[Real Component of $\Psi_l(x,t)$]{
\includegraphics[width=0.4\textwidth]{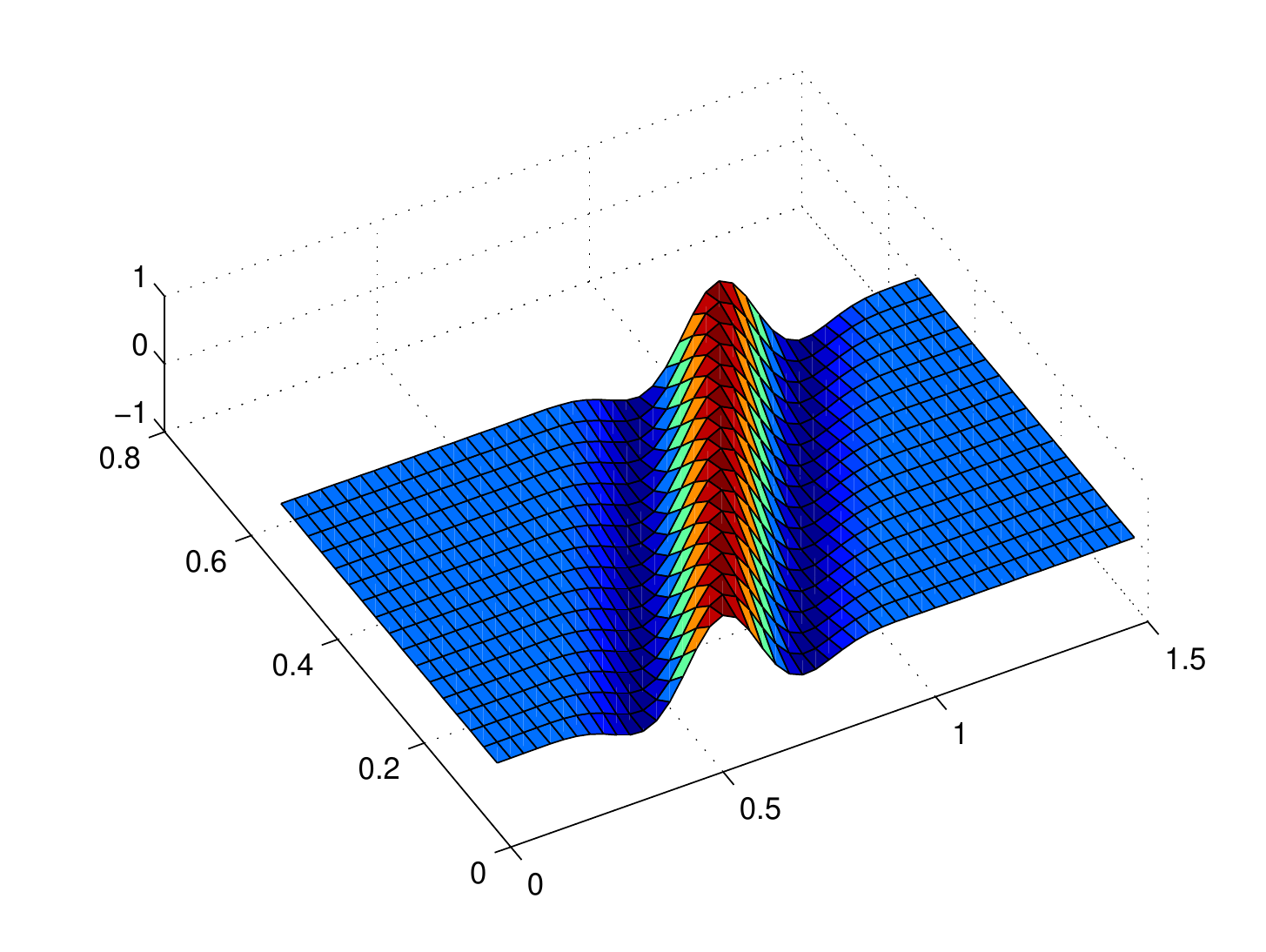} }
\subfigure[Imaginary Component of $\Psi_l(x,t)$]{
\includegraphics[width=0.4\textwidth]{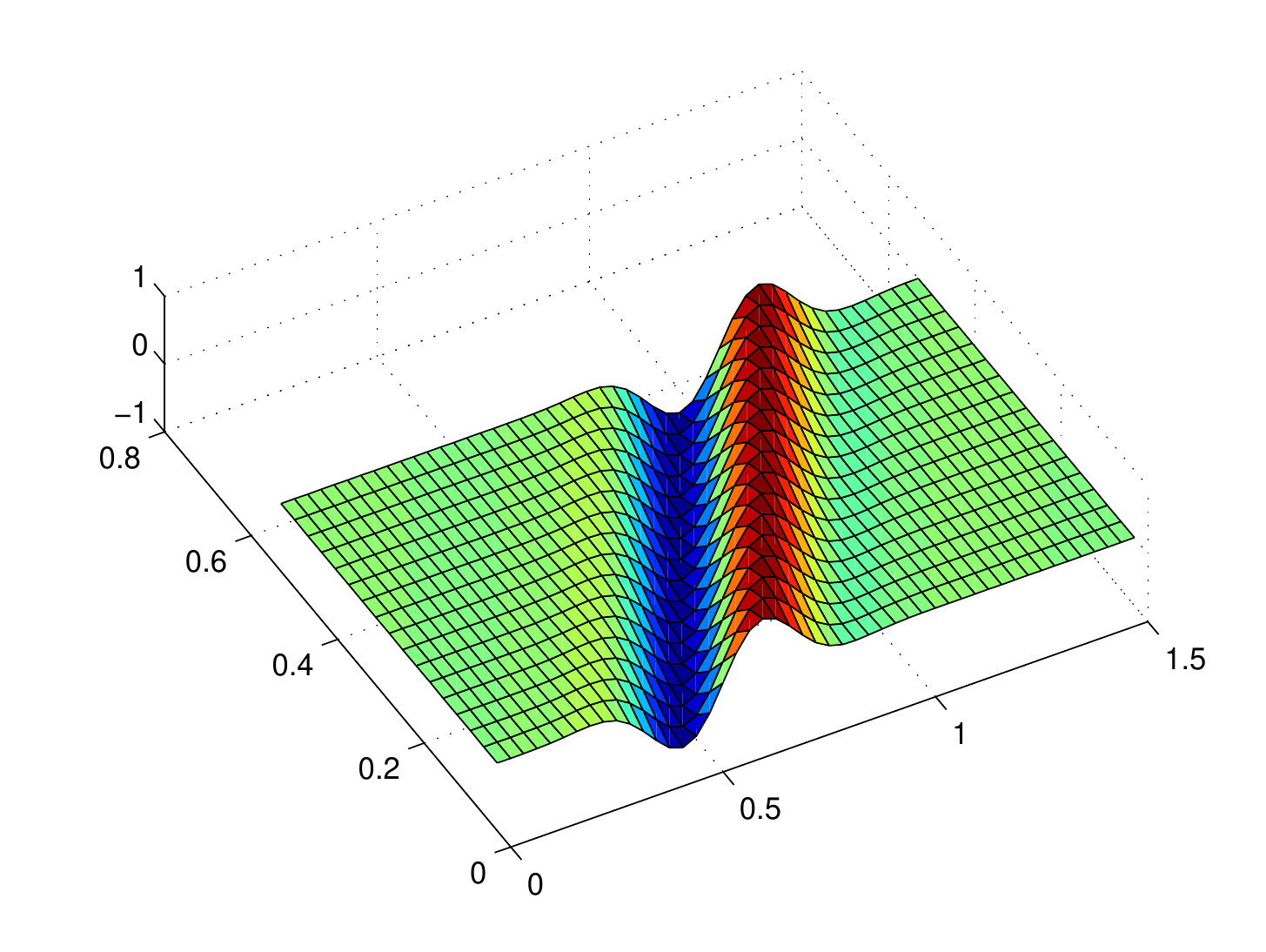} }
\caption{Analytic Solution of the Massless Initial Value Problem}
\label{fig:solution}
\end{figure}

In this section we will show some results from two finite difference approaches and one finite element method to the Dirac equation. These are the central difference method, the staggered finite difference formulation, and the finite element method using triangular finite elements.
In each method we observe significant non-physical effects in the space-time boundary value problem introduced previously.
For the sake of comparison, we refer the reader to figure~\ref{fig:solution} which shows the analytic solution to the space-time boundary value problem proposed above.

\subsection{Central Difference Discretization}

For an implicit implementation of the two-dimensional Dirac equation using the finite difference method, the integral used in the bilinear form above may be replaced with a double summation

\BEA
\sum_i^N \sum_j^N \Phi_{i} \left(i\hbar \BM 1&0 \\ 0&-1 \EM \partial_t  +  i\hbar \BM 0&1 \\ -1&0 \EM \partial_x  -m \right)\Psi_{j}\delta_{i\,j} = 0 \label{eqn:finitedifference} \\
\text{ where } \delta_{i\,j} = 
\left\{ \begin{matrix} 1 & \text{ when } i=j \\
0 & \text{ when } i \neq j
\end{matrix} \right.
\EEA

Often referred to as the naive discretization, the matrix in this approach is built according to equation~\ref{eqn:finitedifference}. Here we use the following central difference definitions for the partial derivative operators.

\BEA
\partial_x \Psi_{x,t} = \frac{1}{2h}\left( \Psi_{x+h,t} - \Psi_{x-h,t}\right) \nonumber \\
\partial_t \Psi_{x,t} =  \frac{1}{2h}\left( \Psi_{x,t+h} - \Psi_{x,t-h}\right) \label{eqn:cdstencil}
\EEA

\begin{figure}[h]
\centering
\subfigure[Real $\Psi_l(x,t)$, $\Delta t =\Delta x$]{
\includegraphics[width=0.4\textwidth]{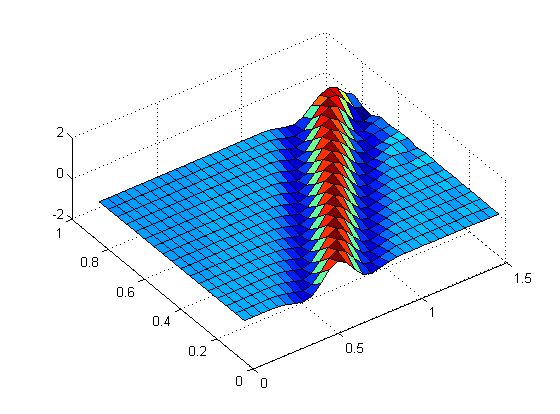} }
\subfigure[Imaginary $\Psi_l(x,t)$, $\Delta t =\Delta x$]{
\includegraphics[width=0.4\textwidth]{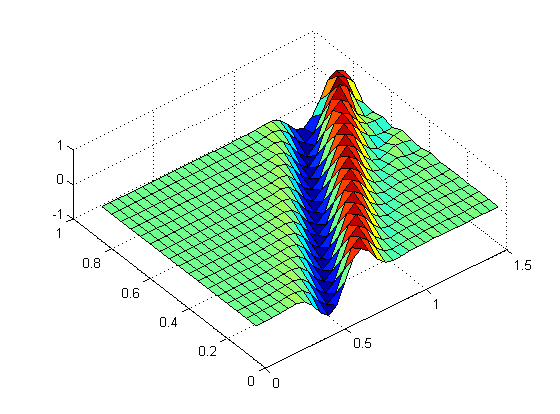} }
\\
\subfigure[Real $\Psi_l(x,t)$, $\Delta t =2\Delta x$]{
\includegraphics[width=0.4\textwidth]{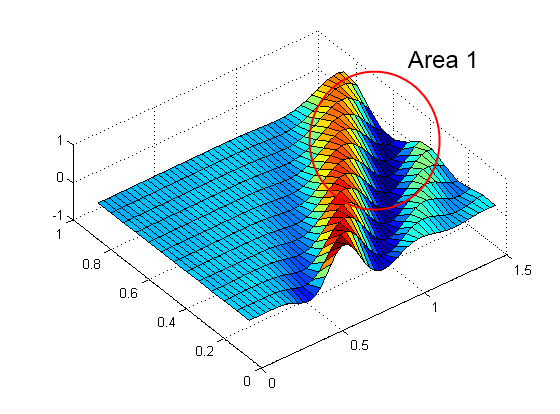} }
\subfigure[Imaginary $\Psi_l(x,t)$, $\Delta t =2\Delta x$]{
\includegraphics[width=0.4\textwidth]{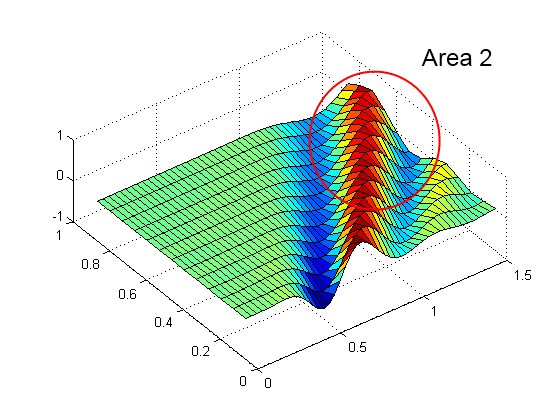} }
\caption{Central Difference Solution. Note that in Area 1 and Area 2, the wave function has shifted rightward, indicating super-luminal propagation, which is forbidden for massless solutions of the Dirac equation}
\label{fig:finitedifference}
\end{figure}

The central finite difference discretization was implemented using the bilinear form of the finite difference method shown in equation~\ref{eqn:finitedifference} and the initial value was introduced via a matrix partitioning scheme.  The result of the balanced difference discretization when applied to the $1+1$ dimensional Dirac initial value problem is shown by Figure \ref{fig:finitedifference}  and Table~\ref{tbl:FDMresult}.

Figure~\ref{fig:finitedifference} shows that the wave function is similar to the analytic solution when the space and time step sizes are equal. However, when unequal step sizes are used the wave shape deteriorates and is shifted rightward, indicating speeds in excess of the speed of light, which is physically impossible.

\begin{table}[ht]
\caption {Numerical Performance of Central Difference Discretization}
\label{tbl:FDMresult}
\begin{center}
\begin{tabular}{|c| c| c| c| c|}\hline
\multicolumn{4}{|c|}{ $\Omega = [0, 1.6]\times[0, 0.8]$ }\\ \hline
\multicolumn{4}{|c|}{$\Delta t =\Delta x$} \\ \hline
Mesh Size			& 	Matrix Size 		&  	BICGSTAB Iterations  	&  Error $\%$ 	\\ 
					\hline
$32\times16$  		&	$1122\times1122$ 	& 	$97.5$			&  $14.38 \%$ 	\\
$48\times24$  		&       	$2450\times2450$ 	& 	$120.5$  			&  $6.28 \%$	\\
$64\times32$  		&	$4290\times4290$ 	& 	$161$				&  $9.04 \%$ 	\\
$80\times40$  		&       	$6642\times6642$ 	& 	$418.5$			&  $15.37  \%$	\\
					\hline \hline
\multicolumn{4}{|c|}{$\Delta t = 2\Delta x$} \\ 
					\hline
$64\times16$  		&       	$2210\times2210$ 	& 	$129$  			&  $76.79 \%$	\\
$80\times20$  		&	$3402\times3402$ 	& 	$250$				&  $56.40 \%$ 	\\
$96\times24$  		&       	$4850\times4850$ 	& 	$307.5$			&  $34.36  \%$	\\
$112\times28$  		&       	$6554\times6554$ 	& 	$578.5$			&  $38.85  \%$	\\
					\hline 
\end{tabular}
\end{center}
\end{table}

Table~\ref{tbl:FDMresult} shows that $L_2$ norm of the error initially initially improves with a finer mesh, but does not improve uniformly and does not appear to tend toward zero with finer mesh spacings. Further, when the spacing is unequal, $\Delta t = 2\Delta x$, the error is substantially larger, which is expected due to its non-physical behavior.

\subsection{Balanced Difference Discretization}

When used to create an explicit propagator, the central difference discretization does not necessarily conserve the probability current of the wave-function. To address this shortcoming the partial derivative stencil in equation~\ref{eqn:cdstencil} is replaced by stecils that are arranged symmetrically with respect to space and time as follows \cite{wessels1999}. 

\BEA
\partial_x \Psi_{x,t} = \frac{1}{4h}\left( \Psi_{x+h,t+h} + \Psi_{x+h,t-h} - \Psi_{x-h,t+h}- \Psi_{x-h,t-h}\right) \nonumber \\
\partial_t \Psi_{x,t} =  \frac{1}{4h}\left( \Psi_{x+h,t+h} + \Psi_{x-h,t+h} - \Psi_{x+h,t-h}- \Psi_{x-h,t-h}\right) \nonumber
\EEA
 
Unlike the original paper  \cite{wessels1999}, where this discretization is used to construct an explicit propagator, our implementation is fully implicit in both time and space.

\begin{figure}[h]
\centering
\subfigure[Real $\Psi_l(x,t)$, $\Delta t =\Delta x$]{
\includegraphics[width=0.4\textwidth]{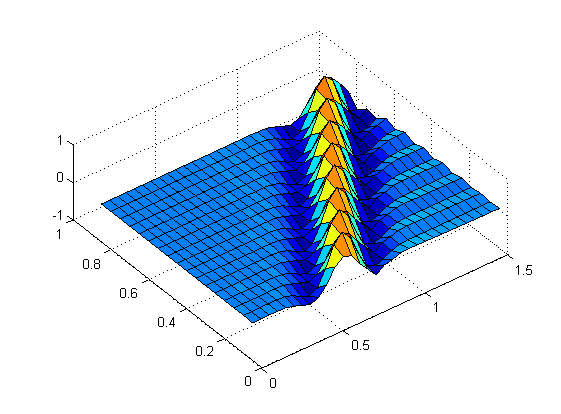} }
\subfigure[Imaginary $\Psi_l(x,t)$, $\Delta t =\Delta x$]{
\includegraphics[width=0.4\textwidth]{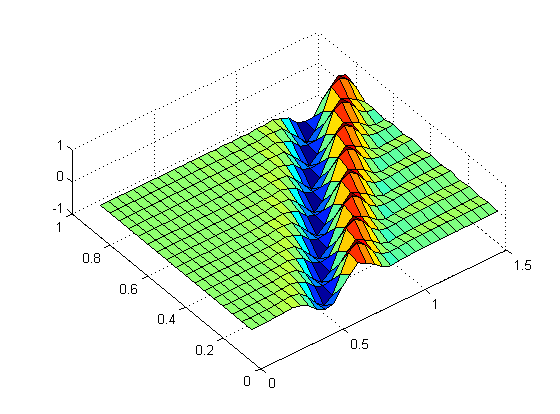} } 
\\
\subfigure[Real $\Psi_l(x,t)$, $\Delta t =2\Delta x$]{
\includegraphics[width=0.4\textwidth]{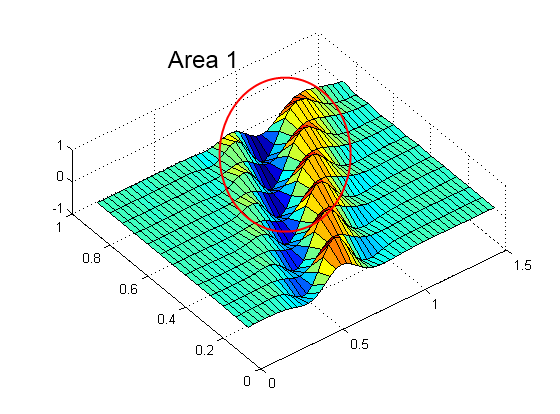} }
\subfigure[Imaginary $\Psi_l(x,t)$, $\Delta t =2\Delta x$]{
\includegraphics[width=0.4\textwidth]{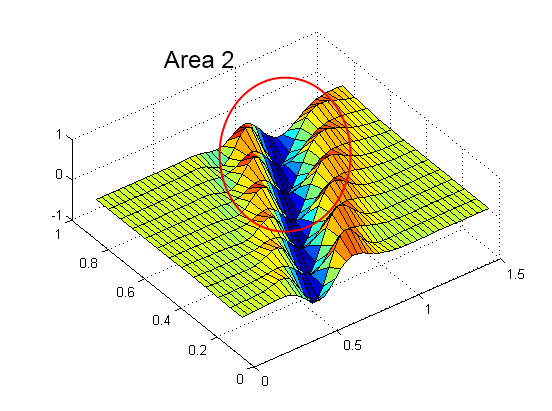} } 
\caption{Balanced Difference Solution. Note that in Area 1 and Area 2, the wavefunction has shifted leftward, indicating sub-luminal propagation, which is inconsistent with the expected behavior of massless solutions}
\label{fig:balanceddifference}
\end{figure}

Figure~\ref{fig:balanceddifference} shows that while the wave shape is choppy, when tested with equal time and space step sizes, it holds a continuous pattern in the overall shape of the analytic solution and the solution shows the correct propagation speed of $v = c$. However, with unequal space and time step sizes the propagation speed is visibly slowed to $v < c$.

\begin{table}[ht]
\caption {Numerical Performance of Balanced Difference Discretization}
\label{tbl:BDMresult}
\begin{center}
\begin{tabular}{|c| c| c| c| c|}\hline
\multicolumn{5}{|c|}{ $\Omega = [0, 1.6]\times[0, 0.8]$ }\\ \hline
\multicolumn{5}{|c|}{$\Delta t =\Delta x$} \\ \hline
Mesh Size		& 	Matrix Size 		&  	BICGSTAB Iterations  	& Residual	&  Error $\%$ 	\\ 
				\hline
$32\times16$  	&	$1122\times1122$ 	& 	$120$			&	.044	&  $27.89 \%$ 	\\
$48\times24$  	&       	$2450\times2450$ 	& 	$863$  		&	.03	&  $19.34 \%$	\\
$64\times32$  	&	$4290\times4290$ 	& 	$53.5$		&	.029	&  $15.30 \%$ 	\\
$80\times40$  	&       	$6642\times6642$ 	& 	$65$			&	.025	&  $11.11  \%$	\\
					\hline\hline
\multicolumn{5}{|c|}{$\Delta t = 2\Delta x$} \\ 
					\hline
$64\times16$  	&       	$2210\times2210$ 	& 	$996$  		&	.13	&  $75.76 \%$	\\
$80\times20$  	&	$3402\times3402$ 	& 	$536$			&	.03	&  $60.29 \%$ 	\\
$96\times24$  	&       	$4850\times4850$ 	& 	$832$			&	.26	&  $45.36  \%$	\\
$112\times28$  	&       	$6554\times6554$ 	& 	$818$			&	.34	&  $38.94  \%$	\\
					\hline 
\end{tabular}
\end{center}
\end{table}

\subsection{Triangular Lagrangian Elements}

\begin{figure}[h]
\centering
\includegraphics[width=0.4\textwidth]{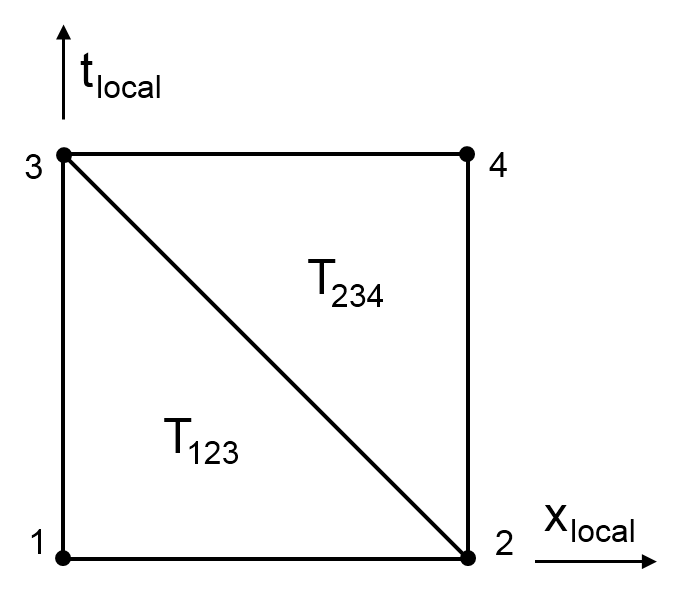} 
\caption{Local Grid Square of a Triangular Finite Element Discretization}
\label{fig:triangles}
\end{figure}

Triangular Lagrangian finite elements represent one of the most commonly used shapes in the finite element method. In this case the nodes of the discrete Dirac equation are arranged on a grid pattern, so each square is composed of two triangular elements as shown in figure~\ref{fig:triangles}. Assuming the single grid element is a unit square with local variables $x$ and $t$, the interpolation polynomials on triangle $T_{123}$ are given by

\BEA
\left.
\begin{matrix}
n_1(x,t) = 1-x-t \nonumber \\
n_2(x,t) = x \nonumber\\
n_3(x,t) = t \nonumber\\
\end{matrix}\right\} &
\text{for } (x,t) \in T_{123} \text{ and} \nonumber\\
\left. \begin{matrix}
n_1(x,t) = 0 \\  
n_2(x,t) = 0 \\
n_3(x,t) = 0 
\end{matrix}\right\}
&
\text{for } (x,t) \text{ elsewhere}\nonumber
\EEA

similarly, interpolation polynomials for $T_{234}$ are given by 

\BEA
\left.
\begin{matrix}
n_2(x,t) = 1-t \nonumber\\
n_3(x,t) = 1 - x \nonumber\\
n_4(x,t) = x + t - 1 \nonumber\\
\end{matrix}\right\} &
\text{for } (x,t) \in T_{234} \text{ and}\nonumber \\
\left. \begin{matrix}
n_2(x,t) = 0 \\  
n_3(x,t) = 0 \\
n_4(x,t) = 0 
\end{matrix}\right\}
&
\text{for } (x,t) \text{ elsewhere}\nonumber
\EEA

Evaluating the finite element the integral in equation~\ref{eqn:weakform} establishes an algebraic relationship between the nodes that for non-boundary elements is equivalent to a finite difference stencil. 
The finite difference stencil for triangular Lagrangian elements are calculated to be

\BEA
\partial_x \rightarrow \frac{1}{36h}
\BM
  -1 & \,1 & \,  \\
-2 & \,0 & \,2 \\
 \, & -1 & \,1
\EM \nonumber & &
\partial_t \rightarrow \frac{1}{36h}
\BM
\,1 & \,2 & \,  \\
-1 & \,0 &  \,1\\
\,  & -2 & -1
\EM
\\
\EEA

The columns of the matrices above correspond to the spatial dimension $x$ and the rows correspond to the temporal dimension $t$. It is apparent that the stencils above are not symmetric with respect to space and time. This means that the choice of element shape may bias the finite difference stencil along the characteristic line $x\, = \,t$ or $x\,=\,-t$, depending on which triangle orientation is chosen. 

\begin{figure}[h]
\centering
\subfigure[Real $\Psi_l(x,t)$, $\Delta t =\Delta x$]{
\includegraphics[width=0.4\textwidth]{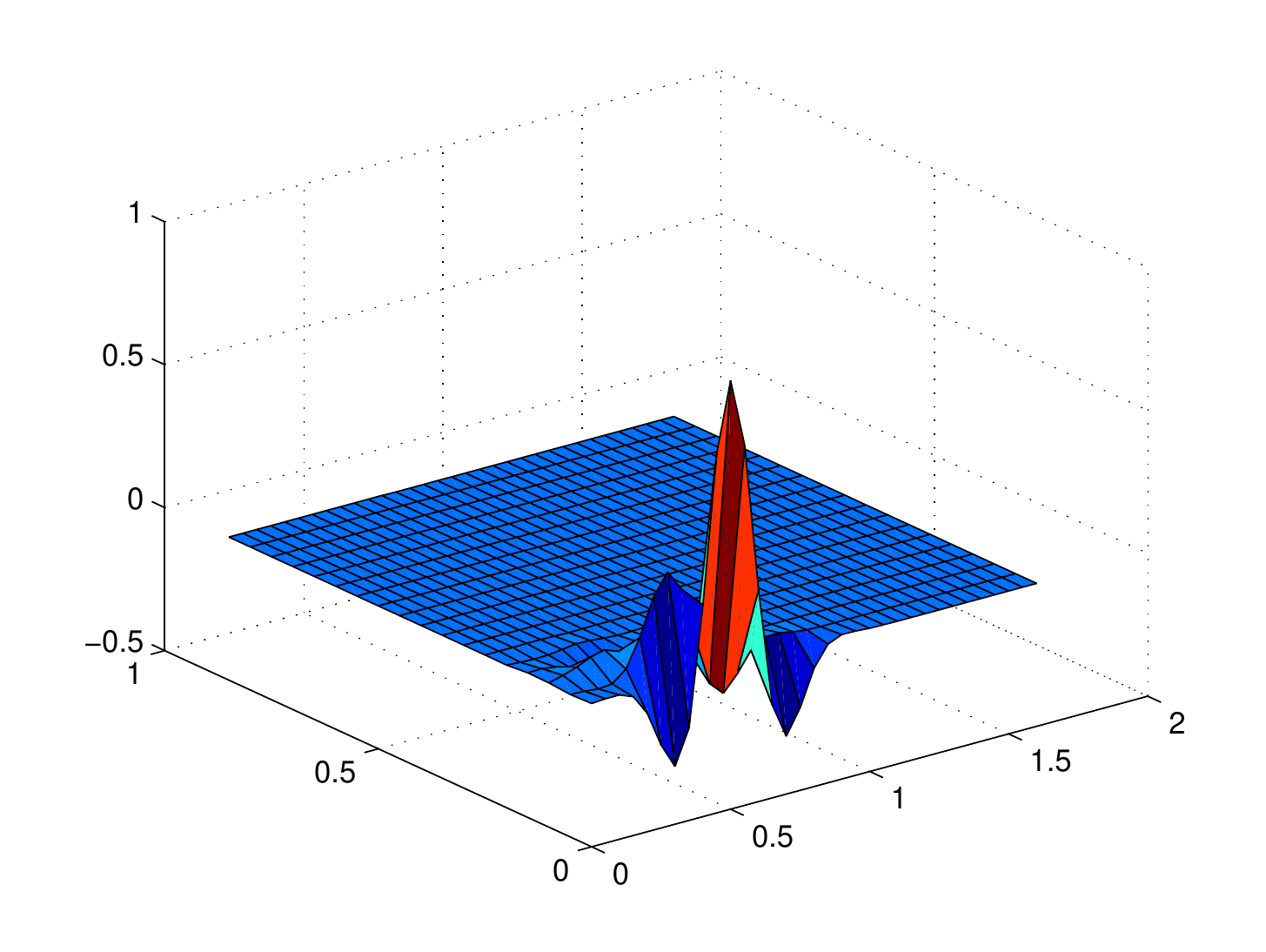} }
\subfigure[Imaginary $\Psi_l(x,t)$, $\Delta t =\Delta x$]{
\includegraphics[width=0.4\textwidth]{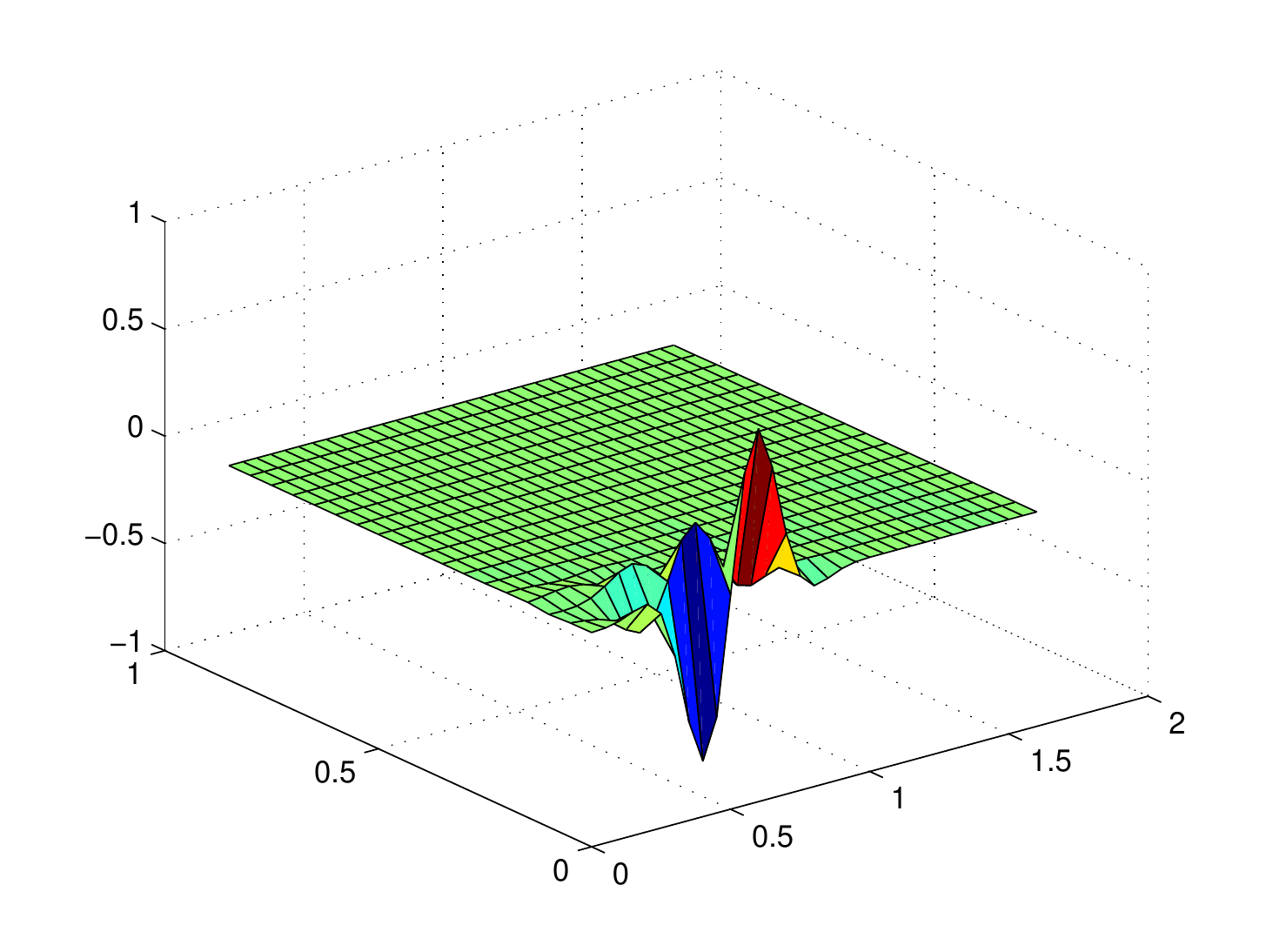} }
\\
\subfigure[Real $\Psi_l(x,t)$, $\Delta t =2\Delta x$]{
\includegraphics[width=0.4\textwidth]{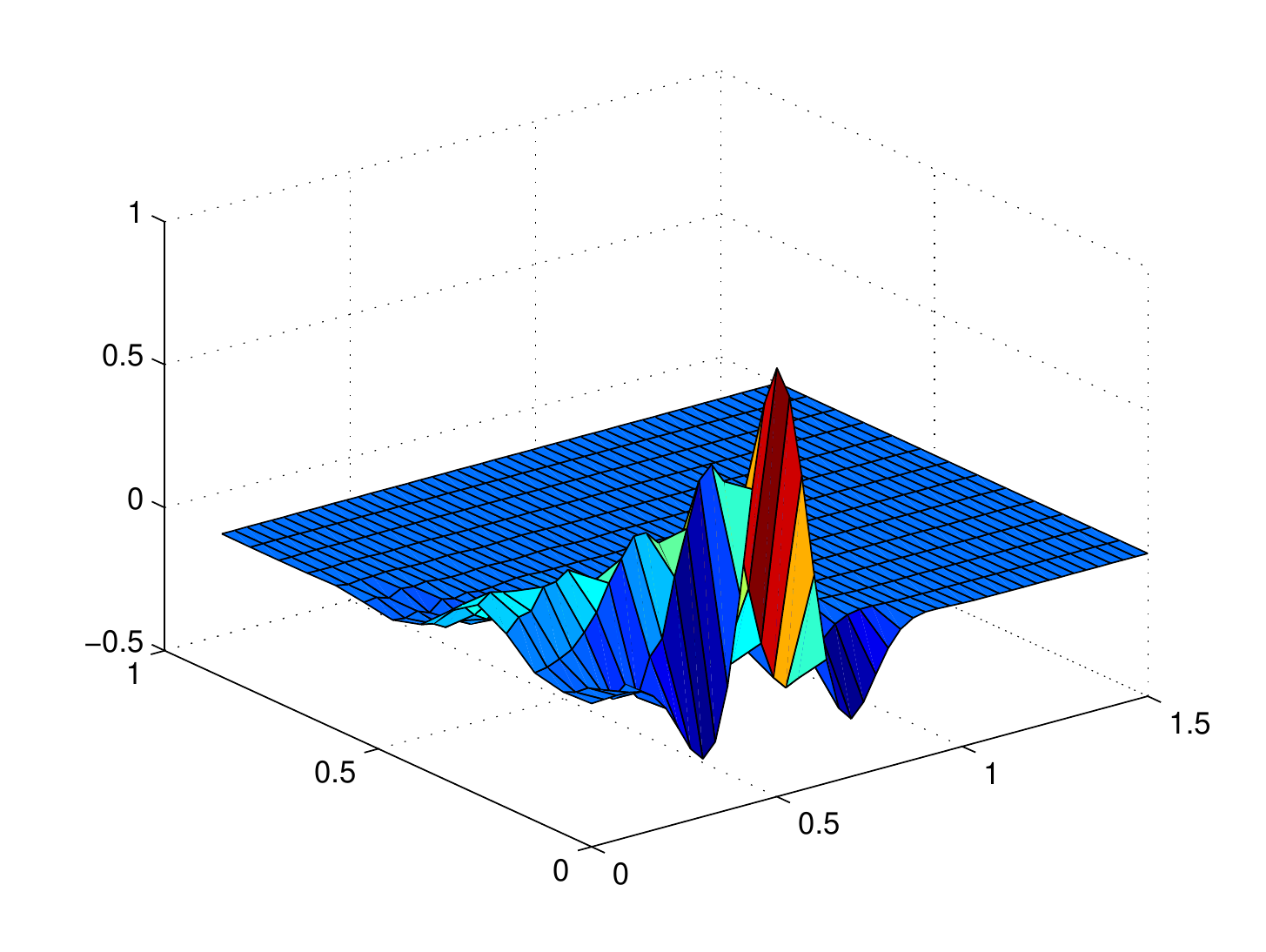} }
\subfigure[Imaginary $\Psi_l(x,t)$, $\Delta t =2\Delta x$]{
\includegraphics[width=0.4\textwidth]{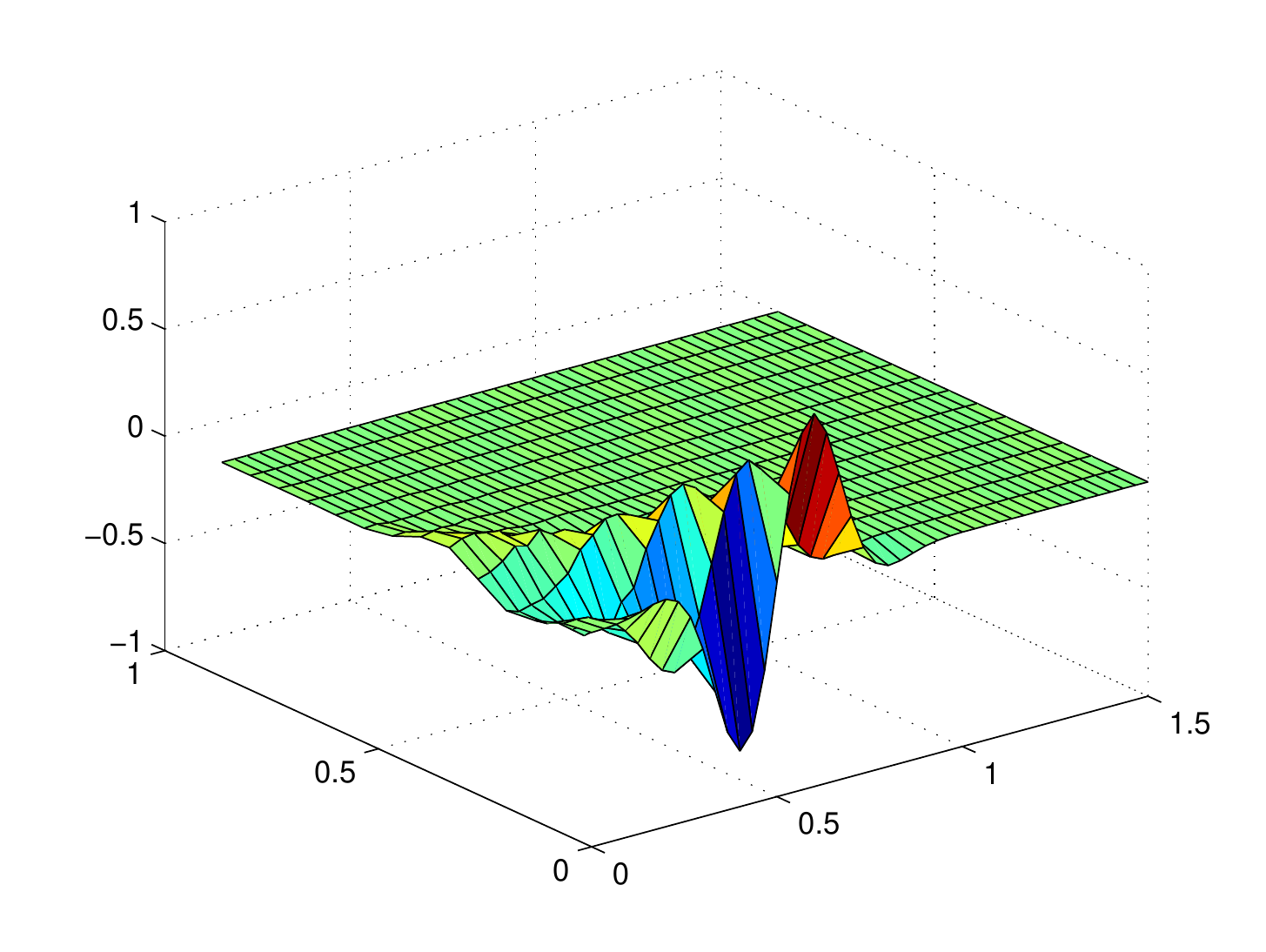} }
\caption{Triangular Lagrangian Element Solution. Note that the calculated wave fades out quickly and propagates in the wrong direction, compared with the analytic solution.}
\label{fig:lagrange}
\end{figure}

Figure~\ref{fig:lagrange} shows that instead of instability, first-order Lagrangian space-time finite elements lose wave amplitude very quickly and appear over-damped. From a physical perspective, the particle is disappearing into space. 
While the precise reason for this disappearance is unknown, it may be due to the finite difference stencils' bias in the opposite direction of particle propagation. When unequal step sizes in space and time were tested, the wave function begins to move to the left in the opposite direction of the analytic solution and at a speed greater than the speed of light.

\begin{table}[ht]
\caption {Numerical Performance of Triangular Lagrangian Elements }
\label{tbl:TLresult}
\begin{center}
\begin{tabular}{|c| c| c| c|}\hline
\multicolumn{4}{|c|}{ $\Omega = [0, 1.6]\times[0, 0.8]$ }\\ \hline
\multicolumn{4}{|c|}{$\Delta t =\Delta x$} \\ \hline
Mesh Size		& 	Matrix Size 		&  	BICGSTAB Iterations  	&  Error $\%$ 	\\ 
					\hline
$32\times16$  	&	$1122\times1122$ 	& 	$2.5$ 				&  $101.6 \%$ 	\\
$48\times24$  	&       	$2450\times2450$ 	& 	$2.5$	  			&  $100.9 \%$	\\
$64\times32$  	&	$4290\times4290$ 	& 	$2.5$				&  $100.6 \%$ 	\\
$80\times40$  	&       	$6642\times6642$ 	& 	$2.5$				&  $100.4  \%$	\\
					\hline \hline
\multicolumn{4}{|c|}{$\Delta t = 2\Delta x$} \\ 
					\hline
$64\times16$  	&       	$2210\times2210$ 	& 	$5$  				&  $108.45 \%$	\\
$80\times20$  	&	$3402\times3402$ 	& 	$7$				&  $108.93 \%$ 	\\
$96\times24$  	&       	$4850\times4850$ 	& 	$4$				&  $106.06  \%$	\\
$112\times28$  	&       	$6554\times6554$ 	& 	$4$				&  $105.49  \%$	\\
					\hline
\end{tabular} 
\end{center}
\end{table}

\section{Tensor Element Based Approaches}
\label{sec:TEBA}

In this section we will show three space-time discretizations that use sqaure shaped "tensor" finite elements along with a  selected basis function space to form the weak Dirac boundary value problem
Each of these approaches shows overall convergence to the shape of the analytic solution without the presence of superluminal, subluminal, or counter-directional wave functions that were present with the previous approaches.

\subsection{Polynomial Hermite Tensor Elements}

These functions are conceptually related to the third-order piecewise Hermite interpolation polynomials given by

\BE
\mathbf{H}(e) = \left\{
\begin{array}{lr}
H_{00}(e) = (1+2e)(1-e)^2 \\
H_{10}(e) = k e(1-e)^2  \\
H_{01}(e) = e^2 (3-2e) \\
H_{11}(e) = k e^2 (e-1) 
\end{array} \right. \nonumber
\EE 

\begin{figure}[h]
\centering
\includegraphics[width=0.4\textwidth]{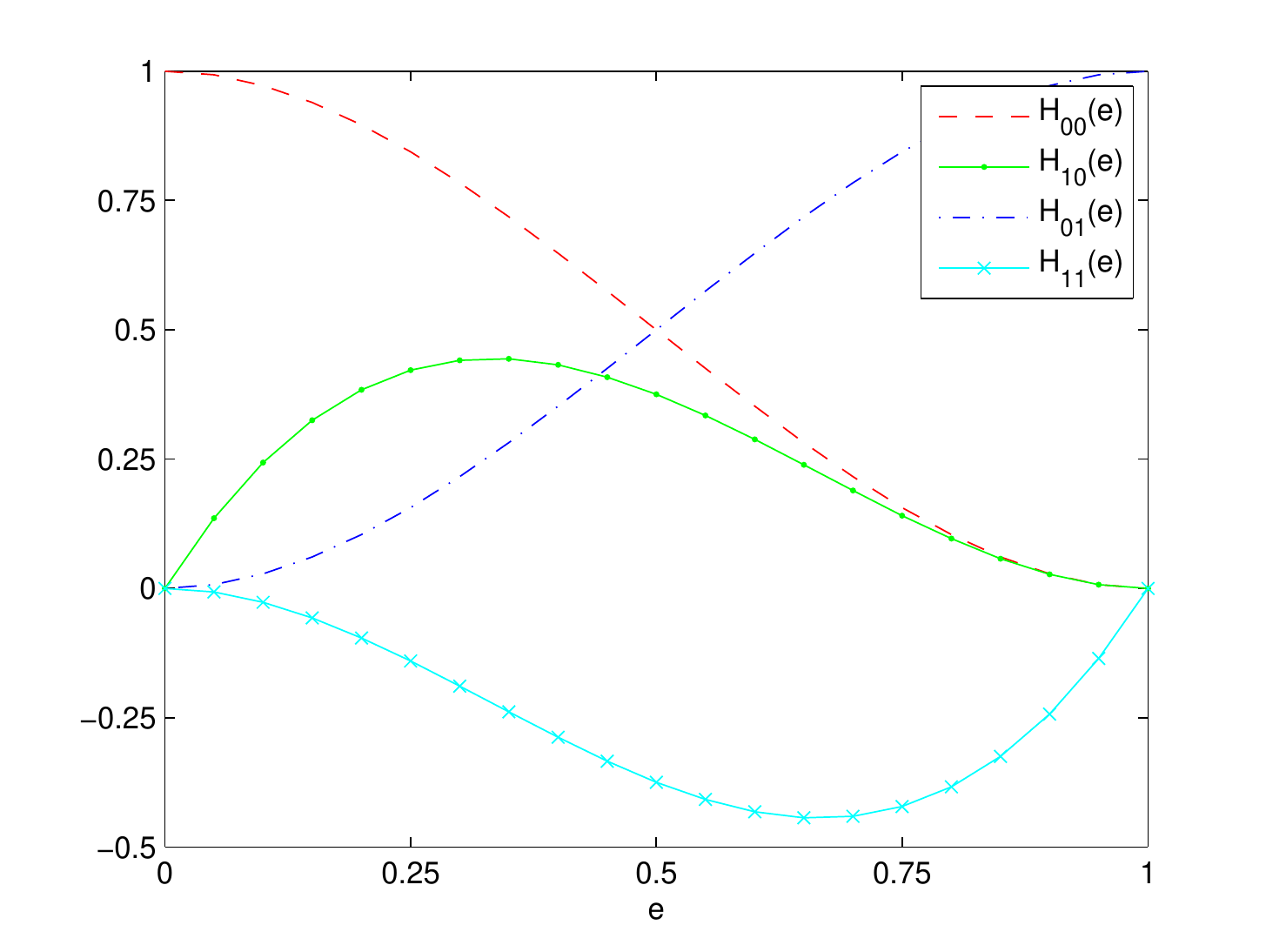}
\caption{$C^1$ Hermite basis functions where $k=3$.}
\label{fig:basisfunctions}
\end{figure}

We use the tensor product $\mathbf H (x) \times \mathbf H (t)$ to create a composite function that has $C^1$ continuity in a two dimensional plane, which is required for analytic solutions in quantum mechanics. 
This set also allows the set of second-order mixed partials to be varied independently. 
However, the continuity of mixed partials has no special physical significance in this case, so these functions are removed in order to reduce the degrees of freedom present in the discrete system.

\begin{figure}[h]
\centering
\subfigure[Real($\Psi_1(x,t)$) Hermite Elements Plot]{
\includegraphics[width=0.4\textwidth]{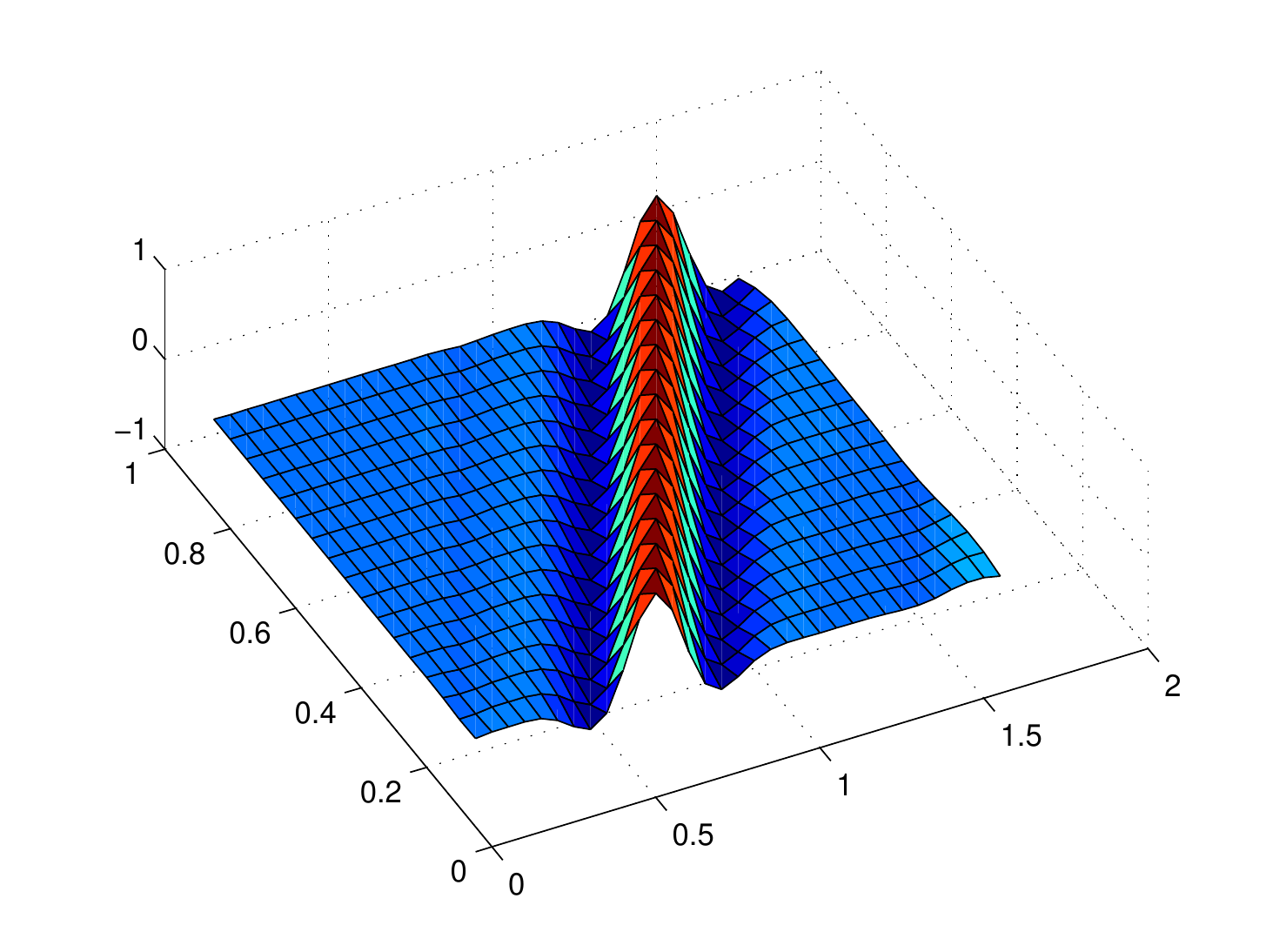} }
\subfigure[Real($\Psi_1(x,t)$) Error Plot]{
\includegraphics[width=0.4\textwidth]{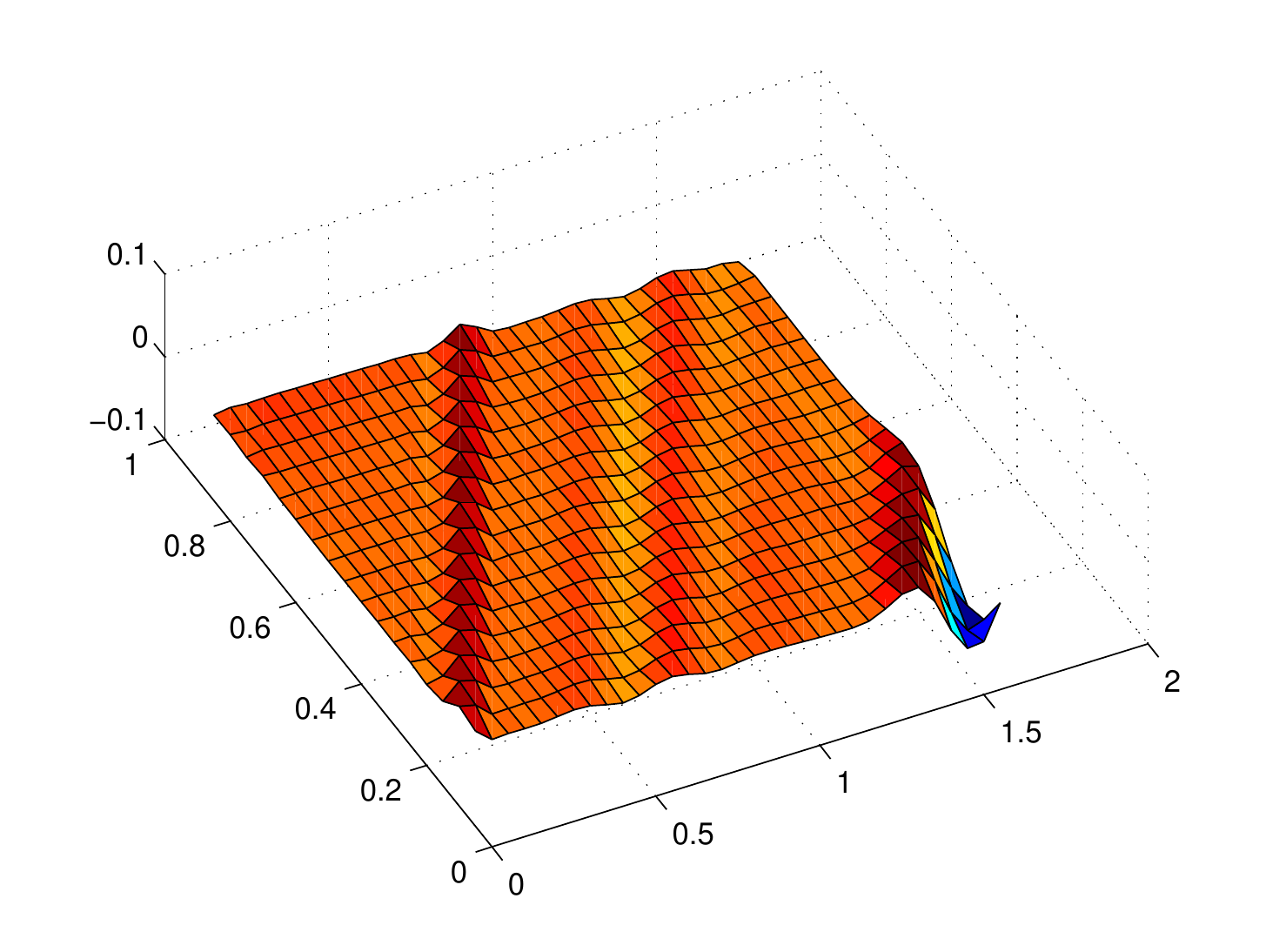} }\\
\subfigure[Real $\Psi_l(x,t)$, Hermite Elements Plot where $\Delta t =2\Delta x$]{
\includegraphics[width=0.4\textwidth]{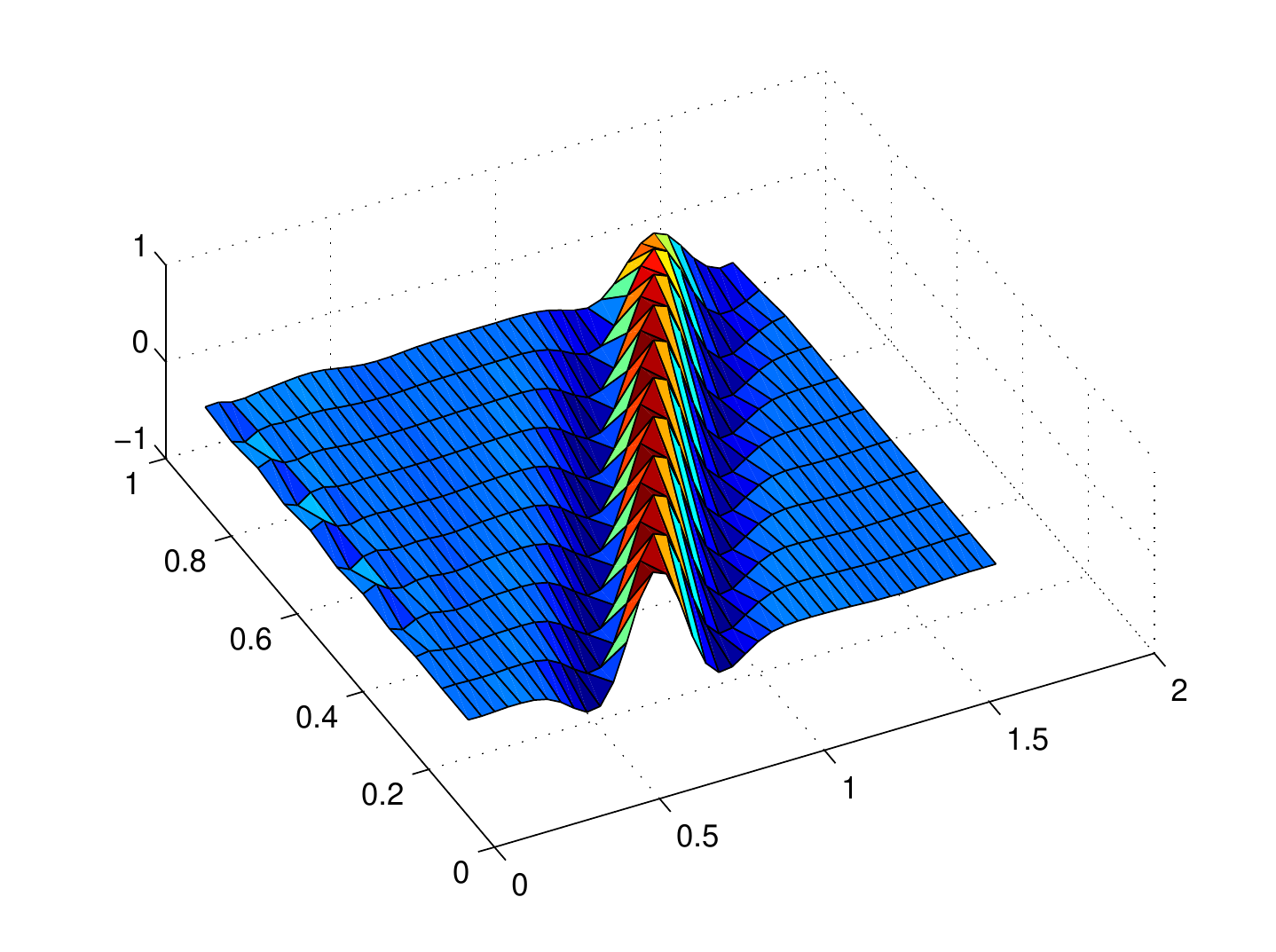} }
\subfigure[Real $\Psi_l(x,t)$, Error Plot where $\Delta t =2\Delta x$]{
\includegraphics[width=0.4\textwidth]{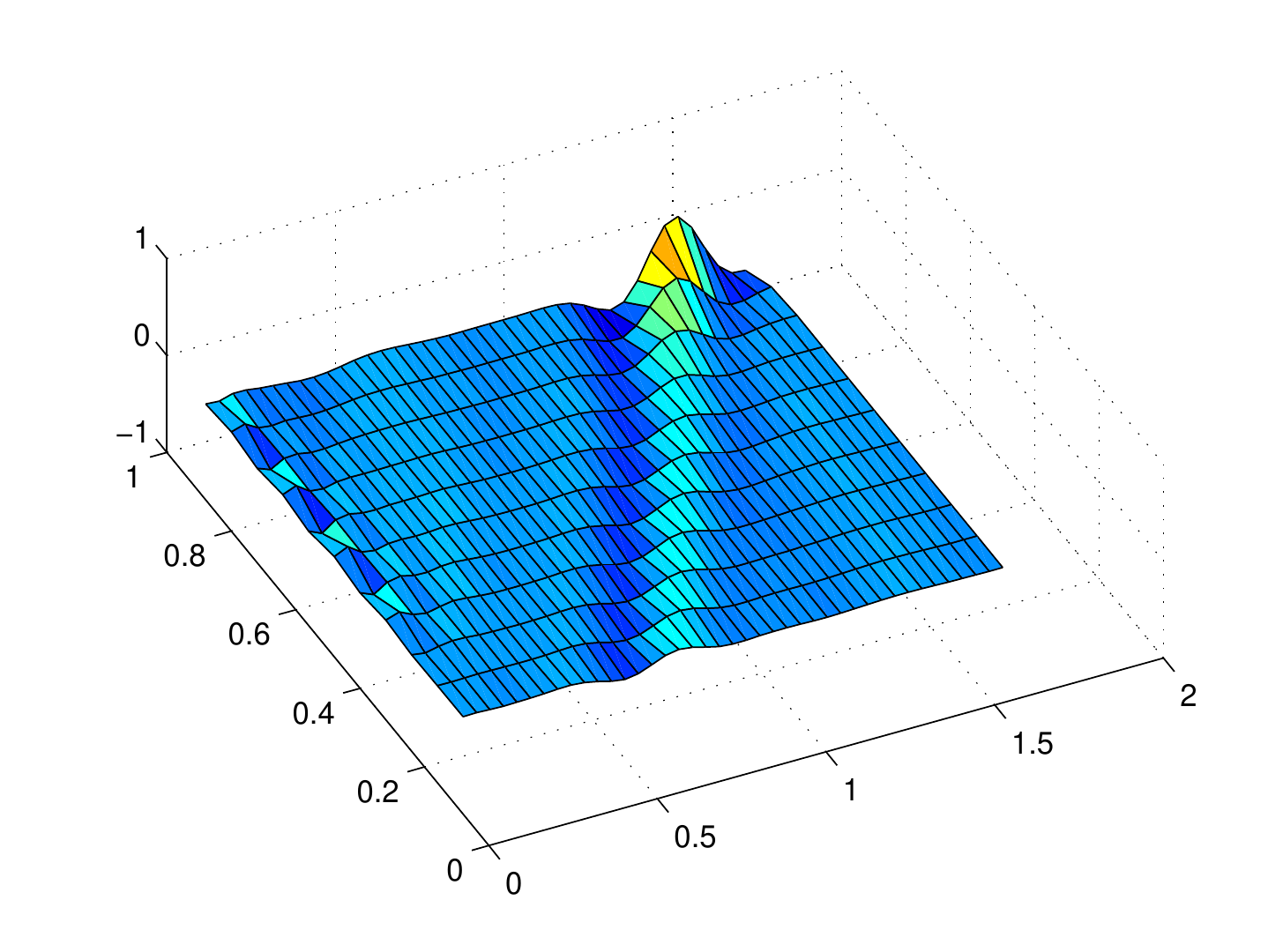} }
\caption{Solution Curve of Polynomial Hermite Tensor Element}
\label{fig:PolyHermite}
\end{figure}

Figure~\ref{fig:PolyHermite} shows that when Hermite tensor elements are used as a function space for the given initial value problem, the overall behavior of the wave-function is consistent with the analytic solution both when $\Delta t = \Delta x$ and when $\Delta t = 2\Delta x$, although the error function is substantial.

\begin{table}[ht]
\caption {Numerical Performance of Hermite Tensor Element Discretization}
\label{tbl:HTEresult}
\begin{center}
\begin{tabular}{|c| c| c| c| c|}\hline
\multicolumn{4}{|c|}{ $\Omega = [0, 1.6]\times[0, 0.8]$ }\\ \hline
\multicolumn{4}{|c|}{$\Delta t =\Delta x$} \\ \hline
Mesh Size			& 	Matrix Size 		&  	BICGSTAB Iterations  	&  Error $\%$ 	\\ 
					\hline
$30\times15$  		&	$2976\times2976$ 	& 	$961$				&  $21.92 \%$ 	\\
$40\times20$  		&       	$5166\times5166$ 	& 	$1801$  			&  $25.64 \%$	\\
$50\times25$  		&	$7956\times7956$ 	& 	$2888$			&  $16.09 \%$ 	\\
$60\times30$  		&       	$11346\times11346$ & 	$4221$			&  $10.97 \%$	\\
$70\times35$  		&       	$15336\times15336$ & 	$5113$			&  $7.35 \%$	\\
$80\times40$  		&       	$19926\times19926$ & 	$6186$			&  $8.81  \%$	\\
					\hline \hline
\end{tabular}
\end{center}
\end{table}

\subsection{Trigonometric Hermite Tensor Elements}

Here we will chose our basis functions for the finite element vectors $\Psi$ to be $C^1$ trigonometric 
functions given by

\BE
\mathbf N(e) = \left\{
\begin{array}{lr}
n_{00}(e) = cos^2(\frac{\pi e}{2} ) \\
n_{10}(e) = k cos(\frac{\pi e}{2})sin(\pi e)  \\
n_{01}(e) = sin^2(\frac{\pi e}{2}) \\
n_{11}(e) = - k sin(\frac{\pi e}{2})sin(\pi e) 
\end{array} \right.\nonumber 
\EE

These are conceptually similar to the Hermite polynomials. As previously, we use the tensor product  $\mathbf N(x) \times \mathbf N(t)$ to create a
 composite function that has $C^1$ continuity.

\begin{figure}[h]
\centering
\subfigure[Real($\Psi_1(x,t)$) Trigonometric Elements Plot]{
\includegraphics[width=0.4\textwidth]{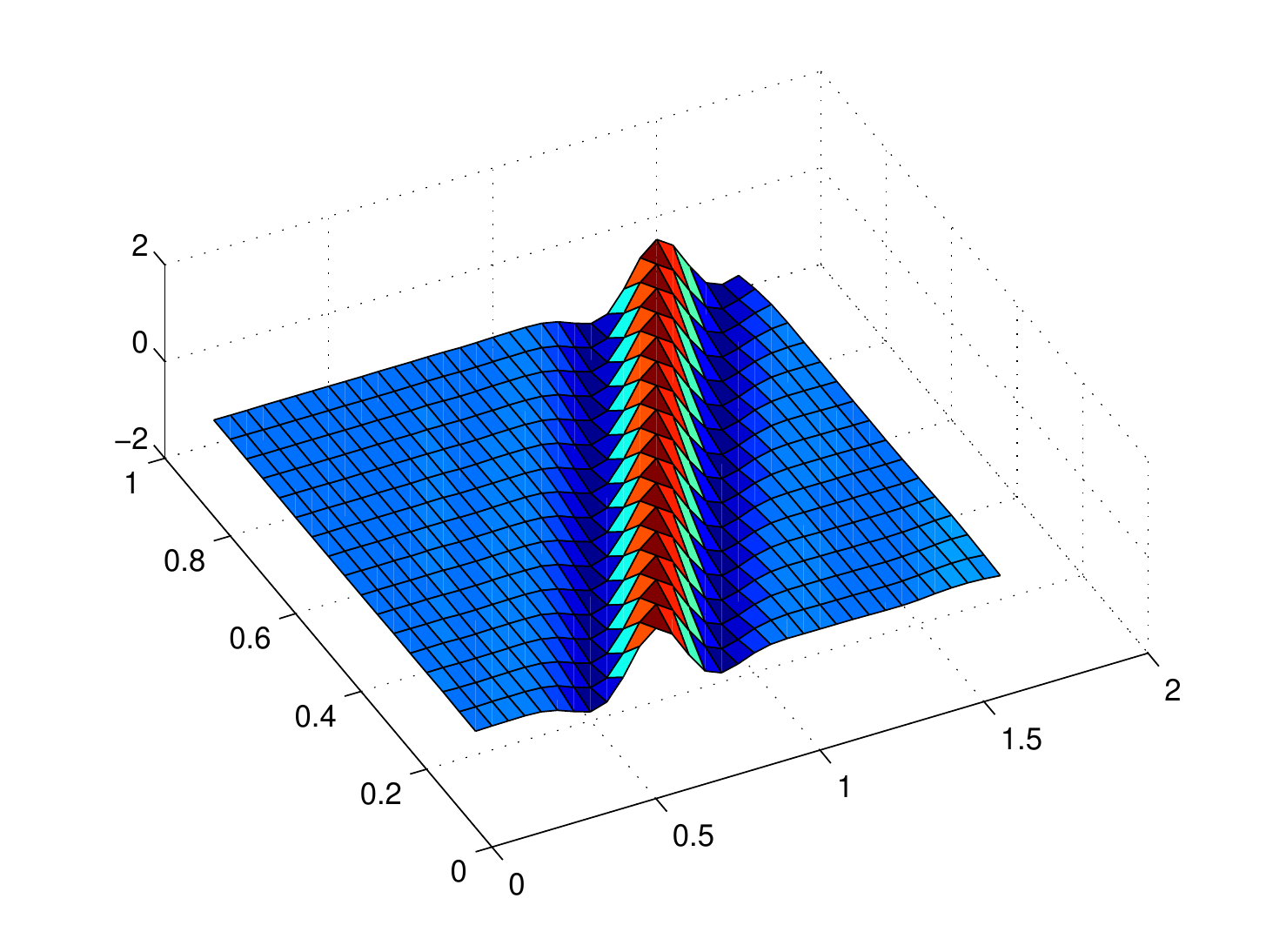} }
\subfigure[Real($\Psi_1(x,t)$) Error Plot]{
\includegraphics[width=0.4\textwidth]{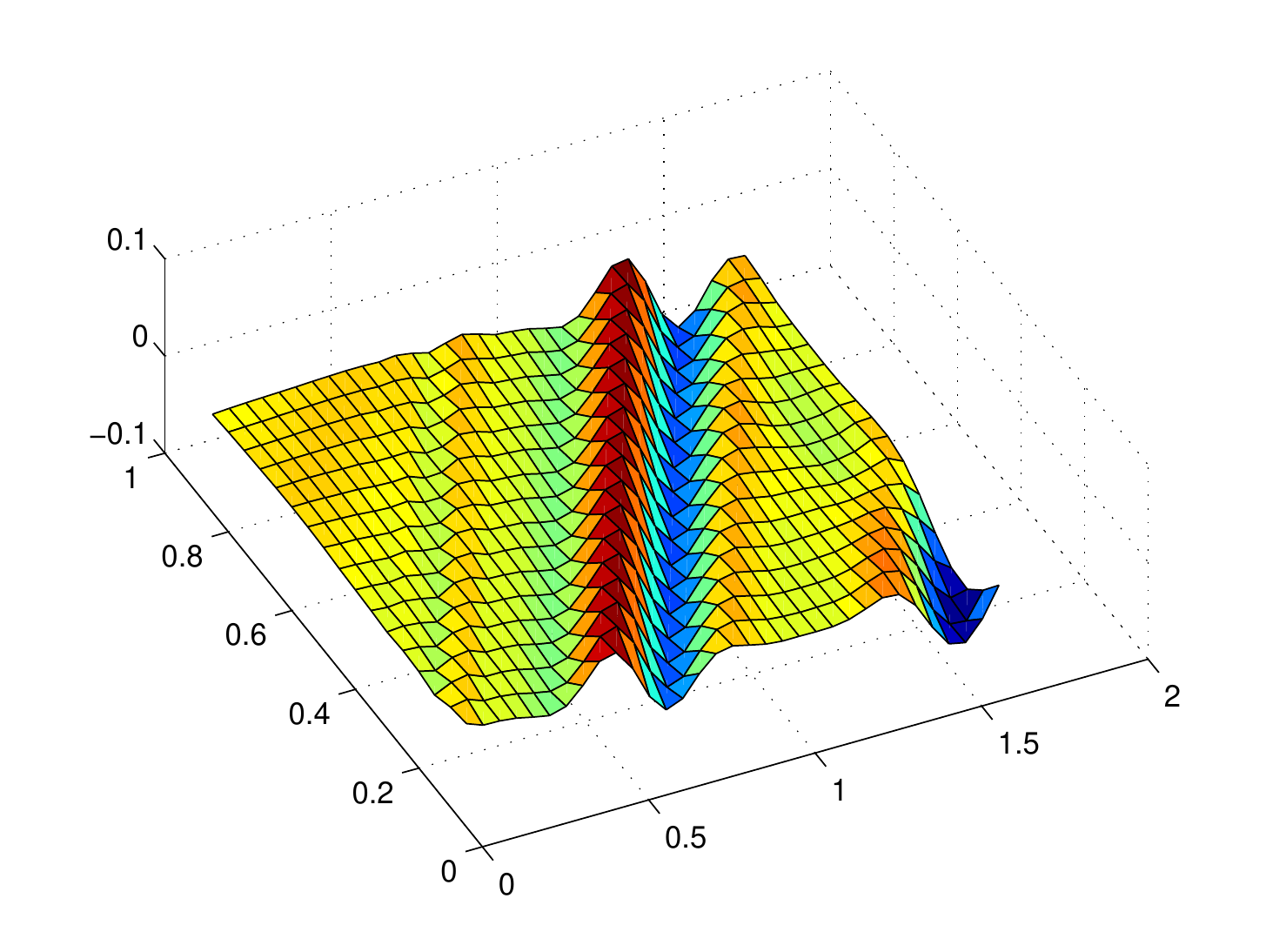} }
\\
\subfigure[Real($\Psi_1(x,t)$) Trigonometric Elements Plot where $\Delta t =2\Delta x$]{
\includegraphics[width=0.4\textwidth]{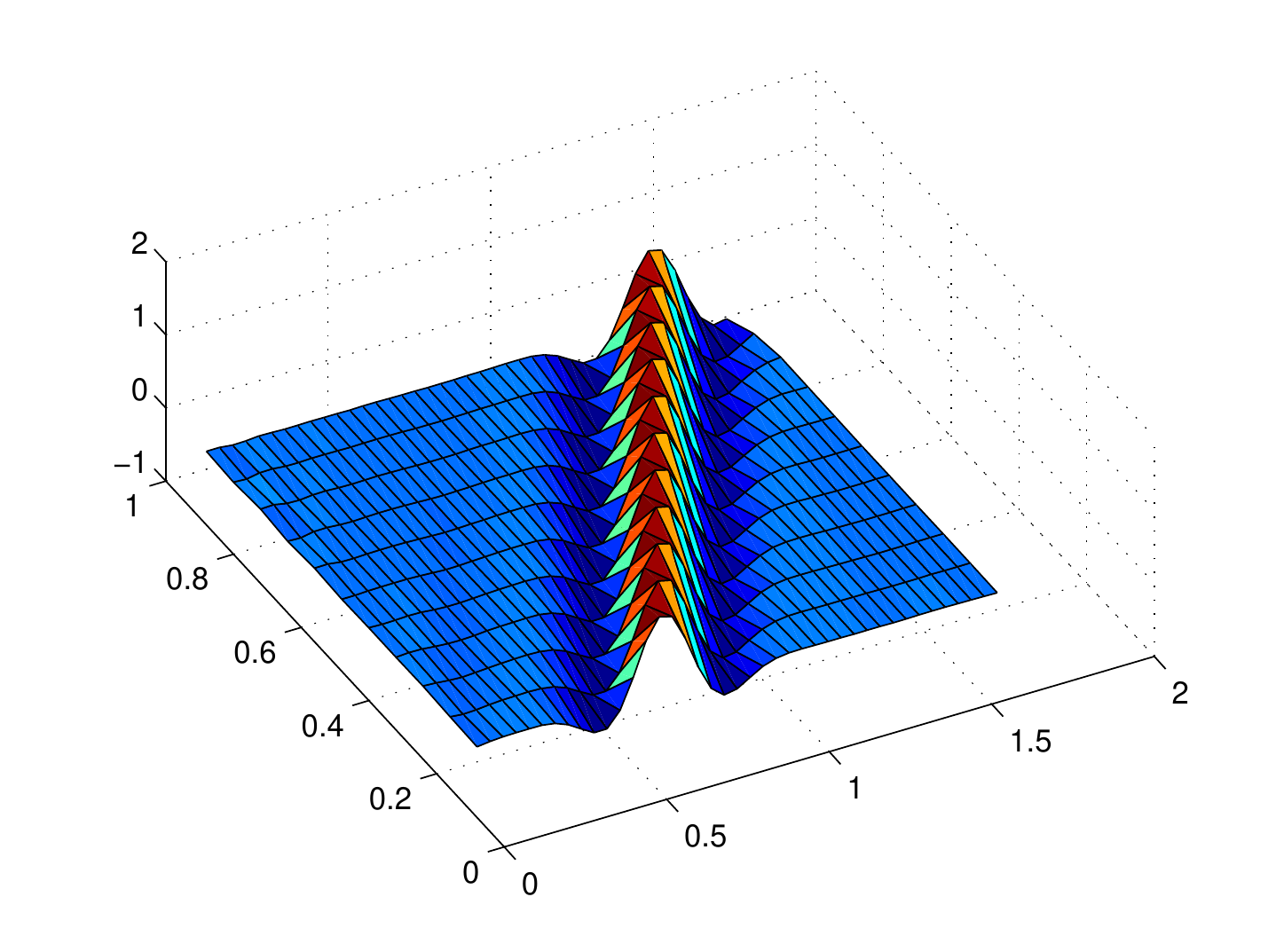} }
\subfigure[Real($\Psi_1(x,t)$) Error Plot where $\Delta t =2\Delta x$]{
\includegraphics[width=0.4\textwidth]{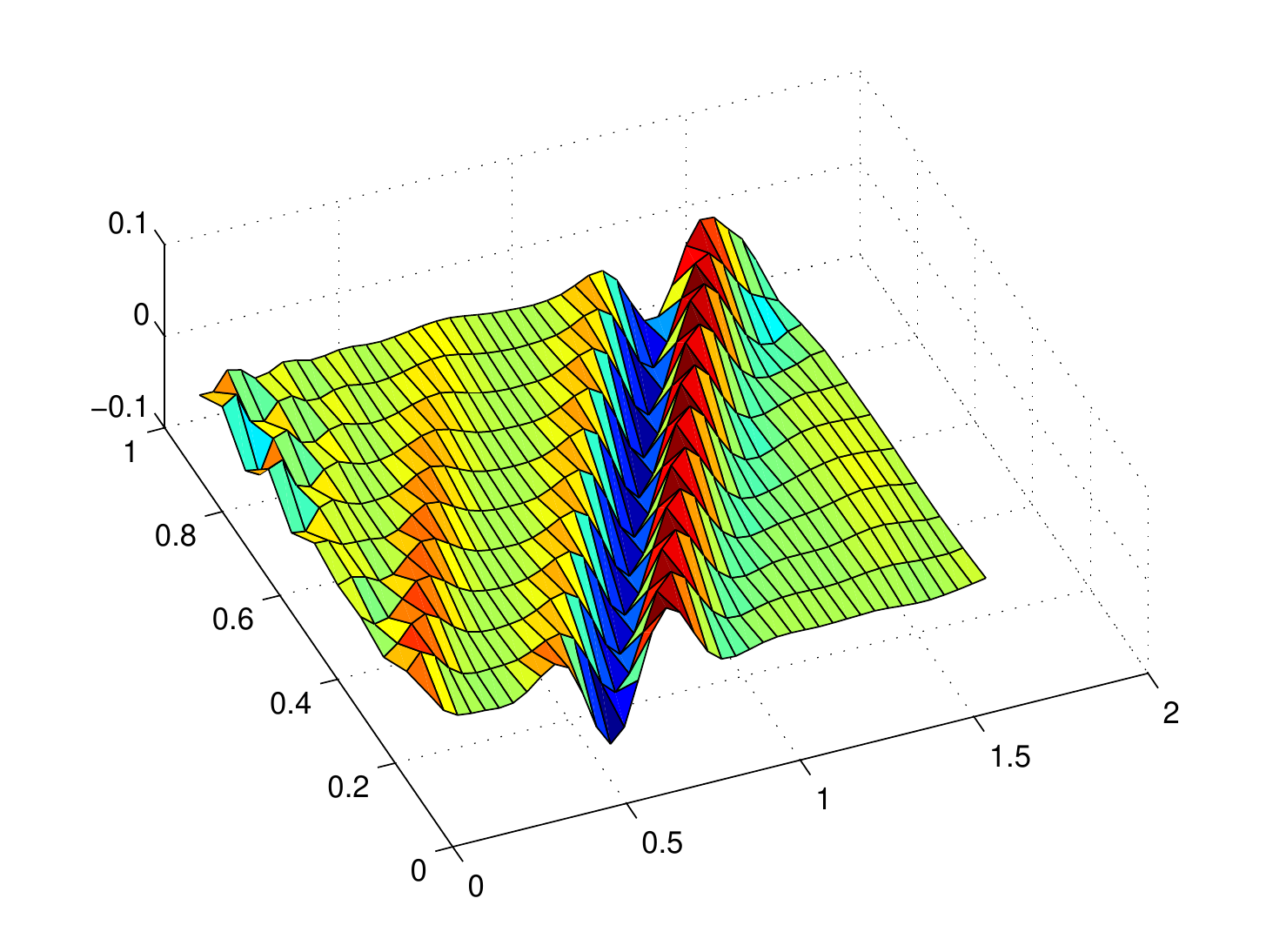} }
\caption{Solution Curve of Trigonometric Tensor Element}
\label{fig:Trig}
\end{figure}

From viewing the solution and error plot in figure~\ref{fig:Trig} we can see that although the overall shape of the solution is very similar to the analytic solution the error wave is relatively large.

\begin{table}[ht]
\caption {Numerical Performance of Trigonometric Tensor Element Discretization}
\label{tbl:TTEresult}
\begin{center}
\begin{tabular}{|c| c| c| c| c|}\hline
\multicolumn{4}{|c|}{ $\Omega = [0, 1.6]\times[0, 0.8]$ }\\ \hline
\multicolumn{4}{|c|}{$\Delta t =\Delta x$} \\ \hline
Mesh Size			& 	Matrix Size 		&  	BICGSTAB Iterations  	&  Error $\%$ 	\\ 
					\hline
$30\times15$  		&	$2976\times2976$ 	& 	$1098$			&  $6.87 \%$ 	\\
$40\times20$  		&       	$5166\times5166$ 	& 	$1992$  			&  $ 5.15\%$	\\
$50\times25$  		&	$7956\times7956$ 	& 	$3021$			&  $ 4.19\%$ 	\\
$60\times30$  		&       	$11346\times11346$ & 	$4375$			&  $ 3.69\%$	\\
$70\times35$  		&       	$15336\times15336$ & 	$5339$			&  $ 3.50\%$	\\
$80\times40$  		&       	$19926\times19926$ & 	$6238$			&  $ 3.18\%$	\\
					\hline \hline
\end{tabular}
\end{center}
\end{table}

\subsection{Linear Lagrangian Tensor Elements}

The function space of the linear Lagrangian elements is formed from the tensor product of the first order Lagrangian interpolation polynomials in the $x$ and $t$ directions. Since they are first order, there is only one degree of freedom per node, or four degrees of freedom per tensor element, making this element type much simpler than the previous elements shown. 

The Lagrangian interpolation polynomials are given by the following expression.
\BE
\begin{array}{lr}
n_1(x,t) = (l-x)(1-t)  \\
n_2(x,t) = x(1-t)\\
n_3(x,t)=(x-1)t\\
n_4(x,t)=x t \nonumber
\end{array}
\EE

\begin{figure}[h]
\centering
\subfigure[Real($\Psi_1(x,t)$) Lagrangian Tensor Elements Plot]{
\includegraphics[width=0.4\textwidth]{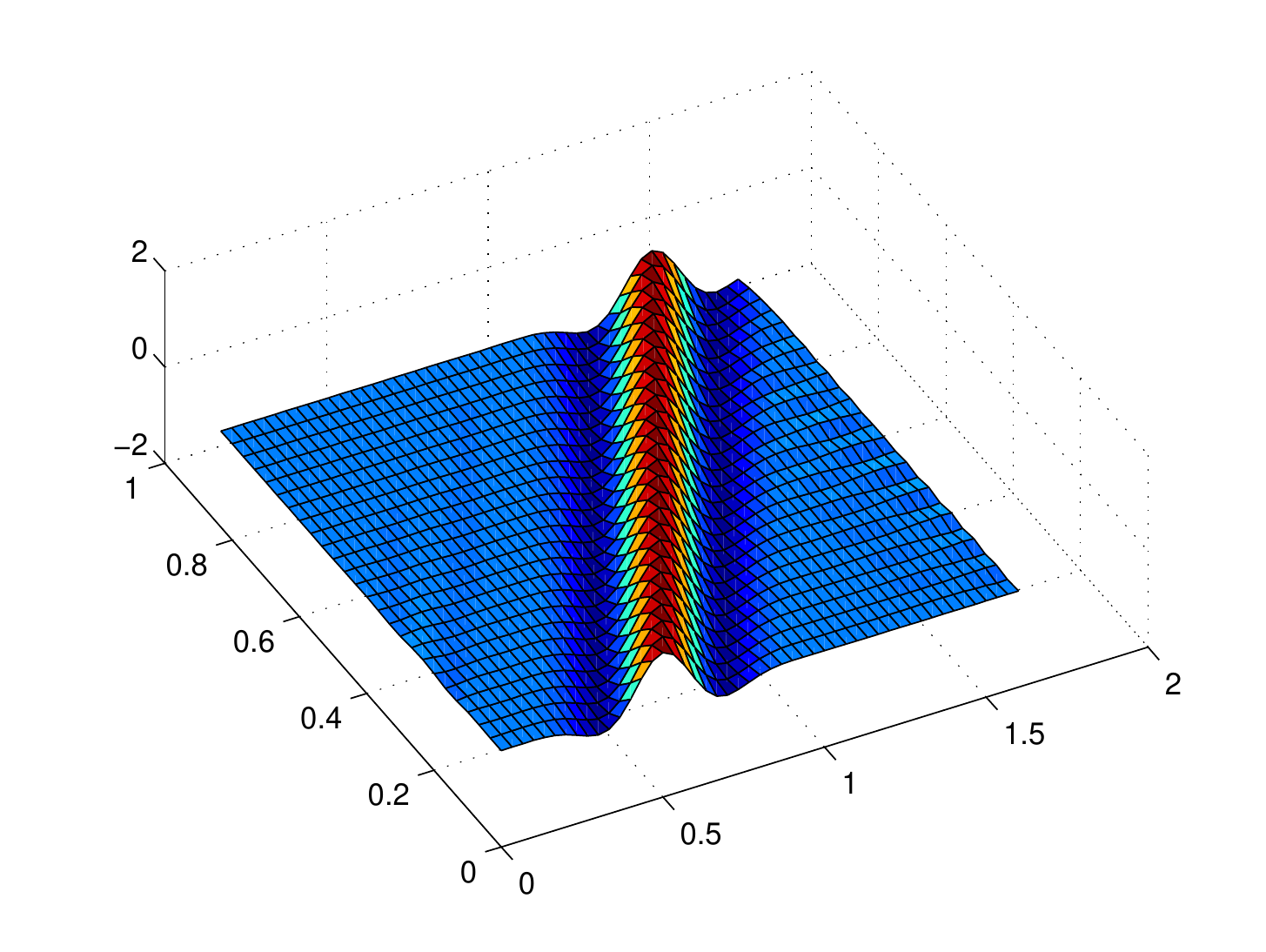} }
\subfigure[Real($\Psi_1(x,t)$) Error Plot]{
\includegraphics[width=0.4\textwidth]{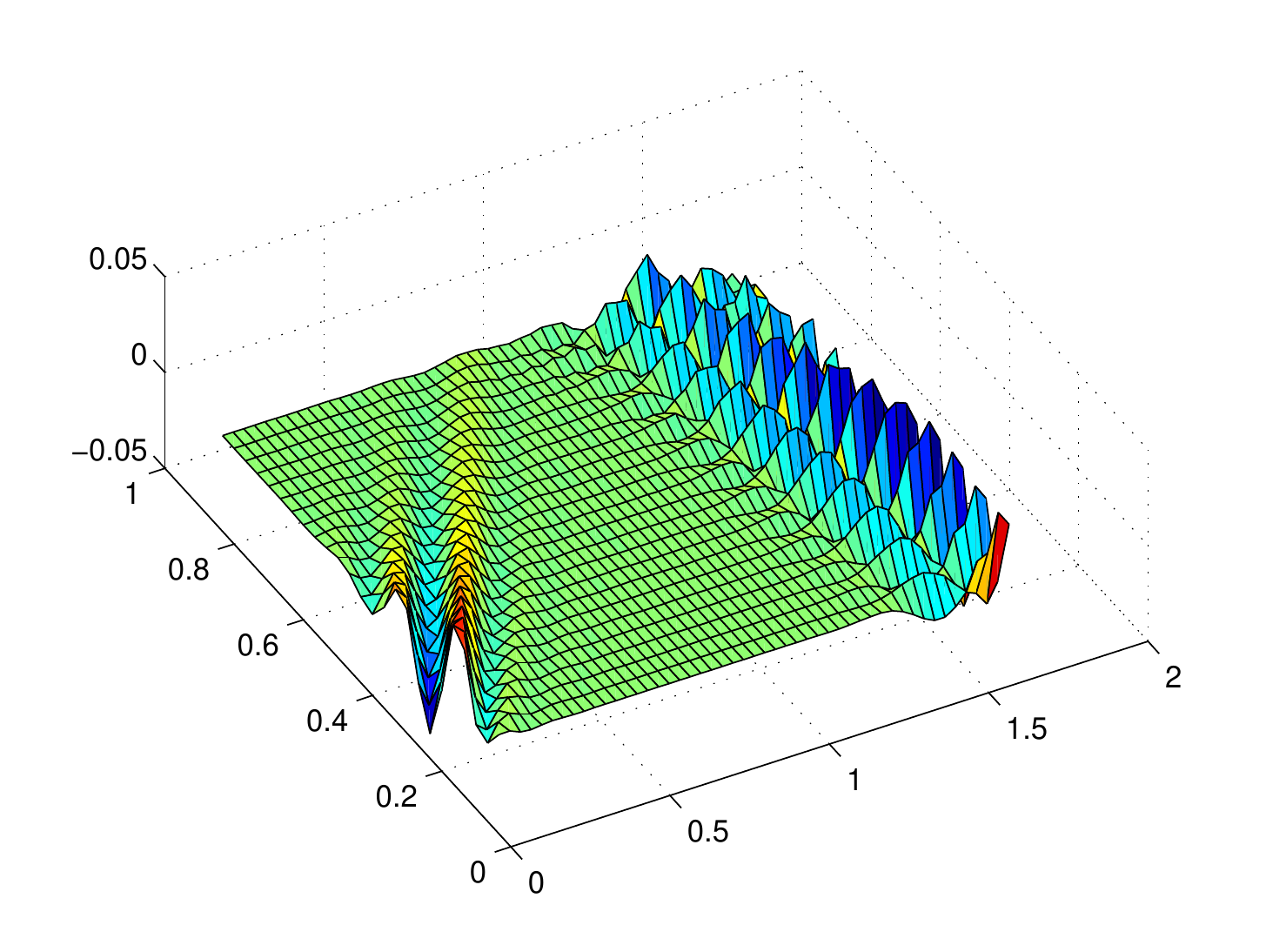} }
\\
\subfigure[Real($\Psi_1(x,t)$) Lagrangian Tensor Elements Plot where $\Delta t =2\Delta x$]{
\includegraphics[width=0.4\textwidth]{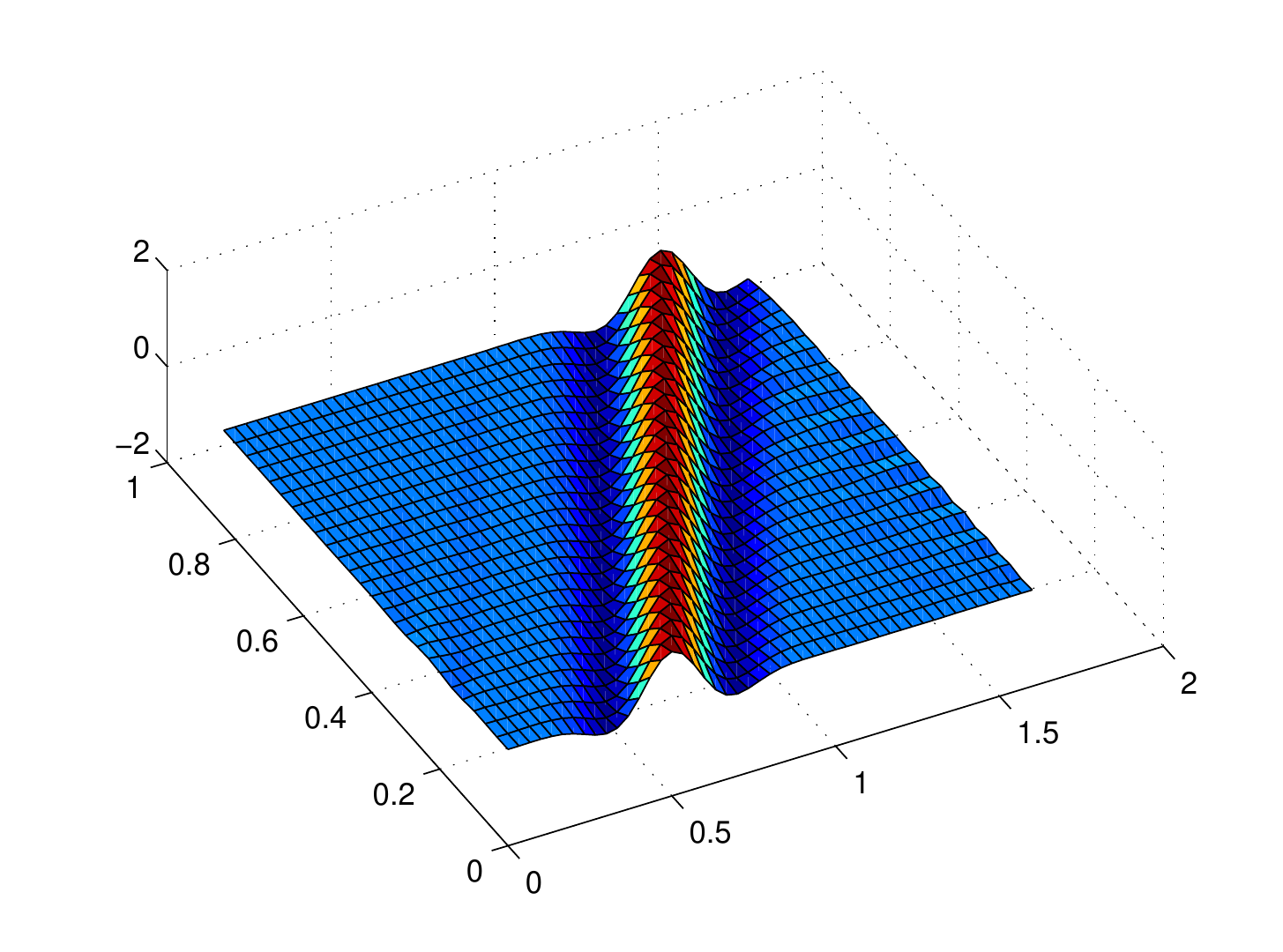} }
\subfigure[Real($\Psi_1(x,t)$) Error Plot where $\Delta t =2\Delta x$]{
\includegraphics[width=0.4\textwidth]{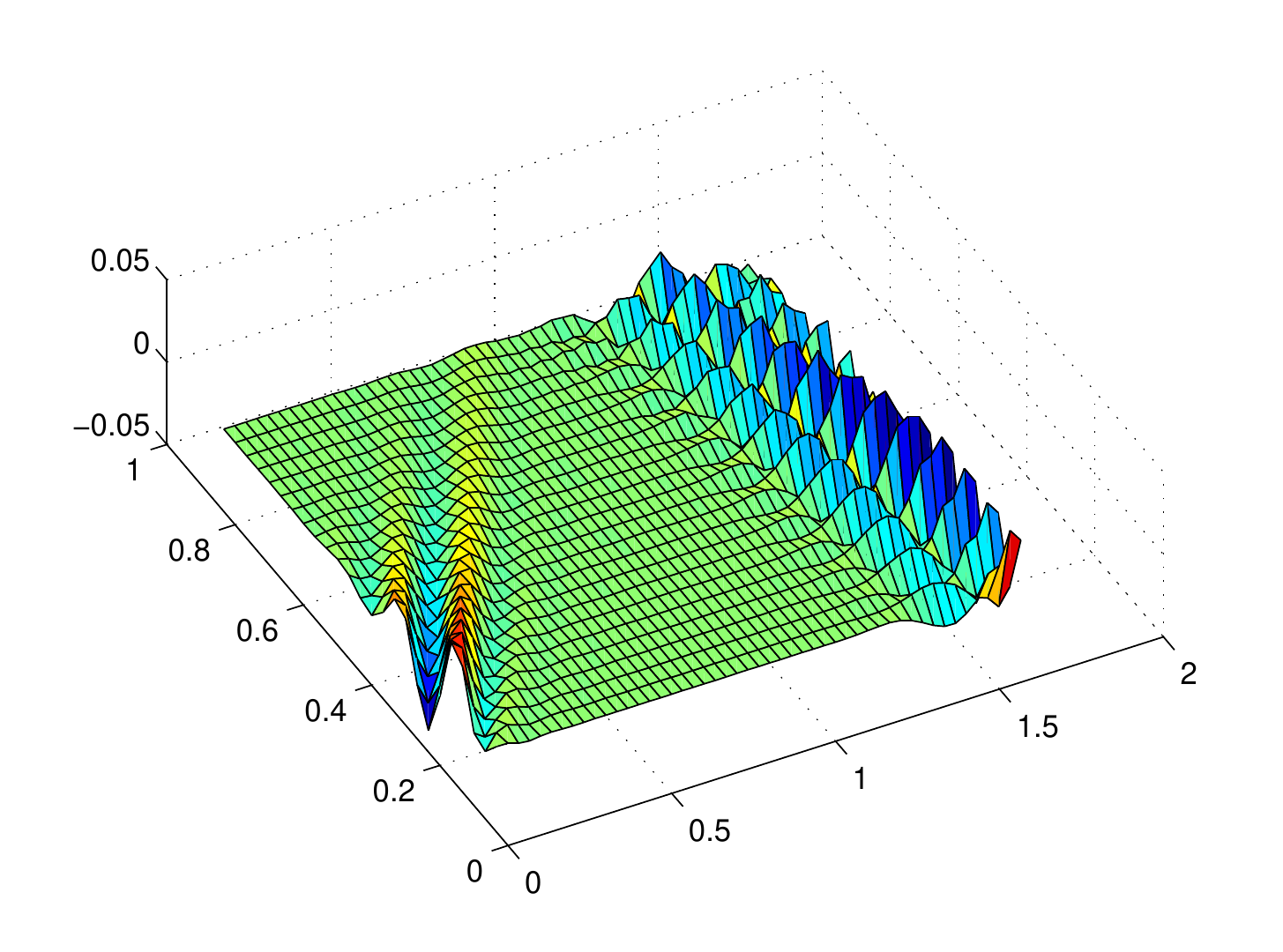} }
\caption{Solution Curve of Lagrangian Tensor Element}
\label{fig:Lagrange}
\end{figure}

From Figure~\ref{fig:Lagrange} we see that the Lagrangian tensor elements also produce a numerical result very close to the analytic solution. This behavior is consistent both when $\Delta t = \Delta x$ and when $\Delta t = 2\Delta x$.

\begin{table}[ht]
\caption {Numerical Performance of Lagrangian Tensor Elements }
\label{tbl:LTresult}
\begin{center}
\begin{tabular}{|c| c| c| c|}\hline
\multicolumn{4}{|c|}{ $\Omega = [0, 1.6]\times[0, 0.8]$ }\\ \hline
\multicolumn{4}{|c|}{$\Delta t =\Delta x$} \\ \hline
Mesh Size		& 	Matrix Size 		&  	BICGSTAB Iterations  	&  Error $\%$ 	\\ 
					\hline
$48\times24$  	&       	$2450\times2450$ 	& 	$331$	  			&  $4.64\%$	\\
$64\times32$  	&	$4290\times4290$ 	& 	$679$				&  $2.86\%$ 	\\
$80\times40$  	&       	$6642\times6642$ 	& 	$899$				&  $2.22\%$	\\
$96\times48$  	&       	$9506\times9506$ 	& 	$1230$  			&  $3.07\%$	\\
$112\times56$  	&	$12882\times12882$ & 	$1569$			&  $3.49\%$ 	\\
$128\times64$  	&       	$16770\times16770$ & 	$2012$			&  $7.67\%$	\\
					\hline \hline
\end{tabular} 
\end{center}
\end{table}

Table~\ref{tbl:LTresult} shows that the Lagrangian tensor elements have substantially lower error, smaller matrix size, greater mesh refinement and more efficient convergence than either the Hermite or trigonometric tensor elements. Unfortunately, the $L_2$ norm of the error actually increases with greater mesh refinement. Possible sources of this remaining error will be analyzed in the following section.

\subsection{Error Analysis of Linear Lagrangian Tensor Elements}

\begin{figure}[h]
\centering
     \includegraphics{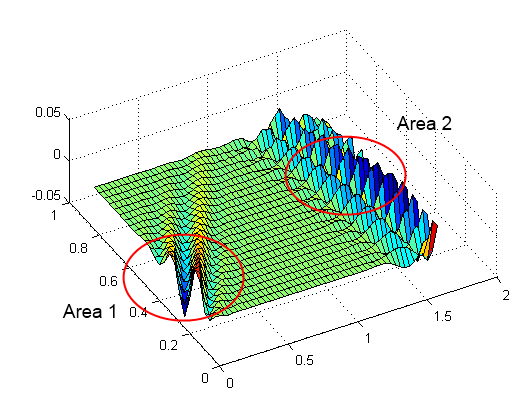} 
     \caption{Error Analysis of Lagrangian Tensor Element. In Area 1 we see boundary error that propagates inward from the natural boundary conditions on the right, left, and rear edges of the domain. In Area 2 the error wave is composed of closely spaced peaks of period 2h, which are effectively invisible to the partial derivative operator for non-boundary nodes}
     \label{fig:LagrangeError}
\end{figure}

From looking at the error wave in figure~\ref{fig:Lagrange} and the analysis in figure~\ref{fig:LagrangeError} we make the following observations about the sources of error.
The first source of error for Lagrangian tensor elements may lie in the stencil of the $\partial_x$ and $\partial_t$ operators. 
For non-boundary nodes, linear Lagrangian tensor elements introduce the following algebraic relationship between nodes. 
\BEA
\partial_x \rightarrow \frac{1}{12h}
\BM
-1 & 0 & 1 \\
-4 & 0 & 4 \\
-1 & 0 & 1
\EM 
& \, &
\partial_t \rightarrow \frac{1}{12h}
\BM
1 & 4 & 1 \\
0 & 0 &  0\\
-1 & -4 & -1
\EM \nonumber
\EEA
It is apparent that the value of the partial derivative as calculated by these operators would approach zero as the period of the wave approaches $2h$.  The error wave across the right-hand side appears to have a period of exactly $2h_t$, making it "invisible" to the discrete form of our partial derivative operator.

A second source of error may come from the boundary conditions across both sides $x=0$ and $x=x_{max}$.
Here an error wave springs \textit{ex-nihilo} from the $x=0$ side and propagates parallel to the solution. 
It may be possible to eliminate such waves by choosing Dirichlet boundary conditions. However, these conditions would imply knowledge solution before the solution is calculated. 
The source of this problem is that the domain sides are not completely contained by the light-cone of our initial condition.

If the domain were sufficiently wide as to preclude the wave from reaching the boundary, it would then be appropriate to apply Dirichlet boundary conditions to the sides of the experiment. However, this would also add siginificant empty space to the domain and computational cost to the experiment.

From the error observations above we draw the following conclusions. One, boundary conditions should utilize light cone causality to ensure a unique solution. Two, momentum and energy operators should be able to detect tightly spaced, erroneous  wave patterns and prevent them from appearing in the solution.. 

\section{Diamond Shapend Tensor Elements}
\label{sec:Diamond}

\begin{figure}[h]
\centering
\subfigure[Diamond Shaped Space-Time Domain]{
	\includegraphics[width=0.4\textwidth]{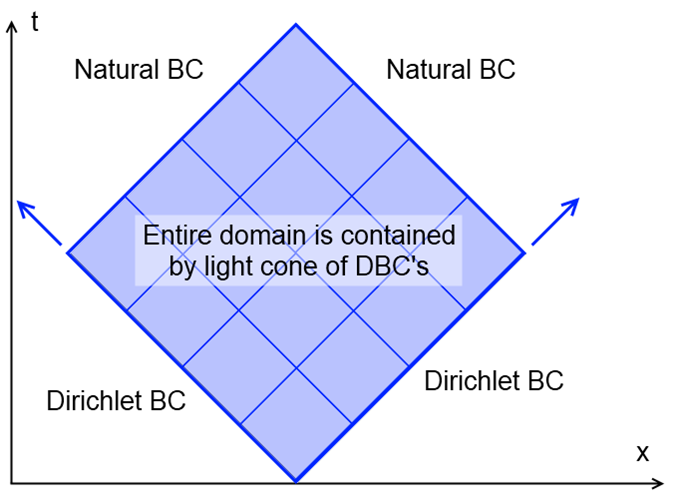} 
	\label{fig:DiamondDomain}}
\subfigure[Diamond Shaped Space-Time Element]{
	\includegraphics[width=0.4\textwidth]{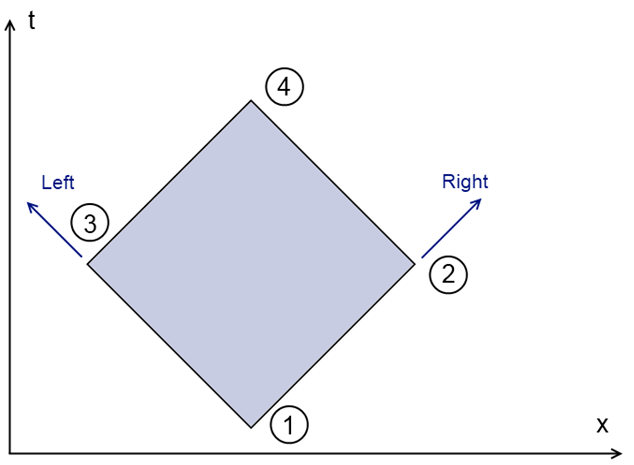} 
	\label{fig:DiamondElement}}
\caption{Diamond Shaped Domain and Single Element Composition}
\label{fig:Diamondcoord}
\end{figure}

In order to reduce the error waves observed in the Lagrangian tensor element discretization, we propose the following element shape for discretizing the 1+1 Dirac equation, shown in Figure~\ref{fig:Diamondcoord}. The $x-t$ plane is rotated $45^\circ$ to create two new axis that will will name "right" and "left". 
This shape gives us two important advantages. 

One, we may impose Dirichlet boundary conditions across both the right and left axis. The entire domain is then contained within the light cone of the "initial" conditions, meaning that the solution should be unique, at least from the physical perspective, since no new information can enter the domain. This is shown graphically in figure~\ref{fig:DiamondDomain}.

Linear interpolation polynomials for the diamond tensor elements may be defined as follows.

\BE
\begin{array}{lr}
n_1(r, l) = (l-r)(1-l)  \\
n_2(r,l) = r(1-l)\\
n_3(r,l)=(r-1)l\\
n_4(r,l)=r l \nonumber
\end{array}
\EE

Two, the partial derivative stencils now become  more complex and should be better able to detect the closely chopped error waves that were present in the Lagrangian tensor element solution. For non-boundary elements, these linear interpolation polynomials introduce the following algebraic relationship between nodes for the two partial derivative operators of the Dirac equation.

\BEA
\partial_x \rightarrow \frac{1}{h}
\BM
 & & 0 & & \\
 & -\frac{1}{3} &  &+\frac{1}{3} & \\
-\frac{1}{6} &  & 0 & &+ \frac{1}{6} \\ 
 & -\frac{1}{3} &  &+\frac{1}{3} & \\
& & 0 & &
\EM
& \,  \, &
\partial_t \rightarrow \frac{1}{h}
\BM
 & &+ \frac{1}{6} & & \\
 &+ \frac{1}{3} &  &+\frac{1}{3} & \\
0 &  & 0 & & 0\\ 
 & -\frac{1}{3} &  &-\frac{1}{3} & \\
& & -\frac{1}{6} & &
\EM \nonumber
\EEA

\begin{figure}[h]
\centering
\subfigure[Real($\Psi_1(x,t)$) Solution Plot]{
	\includegraphics[width=0.4\textwidth]{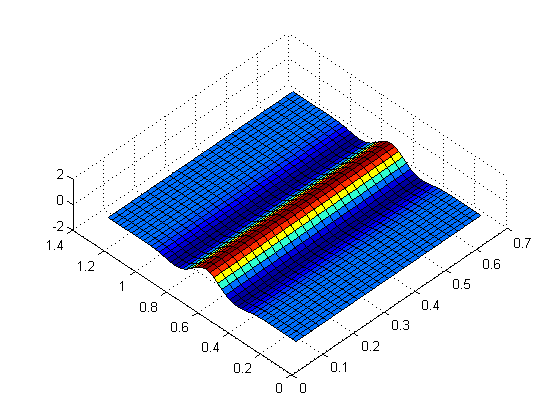} }
\subfigure[Real($\Psi_1(x,t)$) Error Plot]{
	\includegraphics[width=0.4\textwidth]{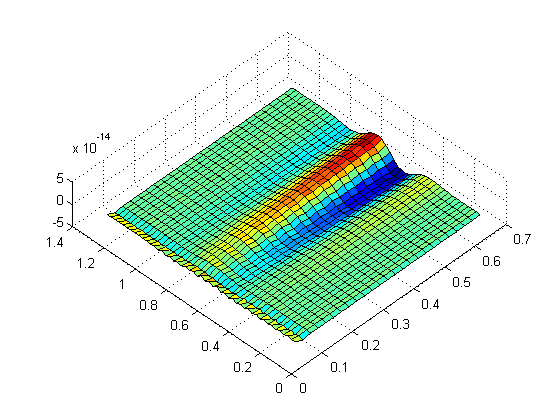} }
\caption{Solution Curve of Diamond Lagrangian Tensor Element. Note the reduction in error scale compared with previous methods tested}
\label{fig:Diamond}
\end{figure}

Using this domain, element shape, and interpolation polynomial set with the weak form defined in equation ~\ref{eqn:weakform} generates the following solution shape shown in Figure~\ref{fig:Diamond}. The magnitude of the error wave is extremely small when compared to the solution, and shows that the finite element solution is nearly exact when one considers the values of the wave-function at the node points. The scale of the error at the nodes is around twelve orders of magnitude lower than the error at the nodes for other methods considered. 

\begin{table}[h]
\caption {Numerical Performance of Diamond Tensor Elements }
\label{tbl:DTresult}
\begin{center}
\begin{tabular}{|c| c| c| c|}\hline
\multicolumn{4}{|c|}{ $\Omega = [0, 1.2]\times[0, 0.8]$ }\\ \hline
\multicolumn{4}{|c|}{$\Delta t =\Delta x$} \\ \hline
Mesh Size		& 	Matrix Size 		&  	BICGSTAB Iterations  	&  Error $\%$ \\ 
					\hline
$24\times48$  	&       	$2450\times2450$ 	& 	$305$	  			&  $1.28\%$	\\
$32\times64$  	&	$4290\times4290$ 	& 	$621$				&  $0.71\%$ \\
$40\times80$  	&       	$6642\times6642$ 	& 	$883$				&  $0.46\%$	\\
$48\times96$  	&       	$9506\times9506$ 	& 	$1177$  			&  $0.32\%$	\\
$56\times112$  	&	$12882\times12882$ & 	$1593$			&  $0.23\%$ \\
$64\times128$  	&       	$16770\times16770$ & 	$1897$			&  $0.18\%$	\\
					\hline 
\end{tabular} 
\end{center}
\end{table}

We note that the solution shows substantially lower error when compared to any of the methods previously presented. This is further confirmed by our results in ~\ref{tbl:DTresult}. Here we see that the $L_2$ norm of the error is much lower than for the other methods tested, and that the numerical simulation converged more quickly as well. In the case of a the $64\times128$ element matrix the diamond tensor element solution converged in $1897$ iterations versus $2012$ iterations for the Lagrangian tensor elements, and the $L_2$ norm of the error was $0.18\%$ (Table ~\ref{tbl:DTresult}) for the diamond tensor elements versus $7.67\%$ (Table ~\ref{tbl:LTresult}) for the Lagrangian tensor elements.

Finally, as with the other space-time tensor element approaches, no superluminal or subluminal behavior was observed when tested with unequal space and time spacings.

\subsection{Rotation Tests}

\begin{figure}[h]
\centering
\includegraphics[width = 0.3\textwidth]{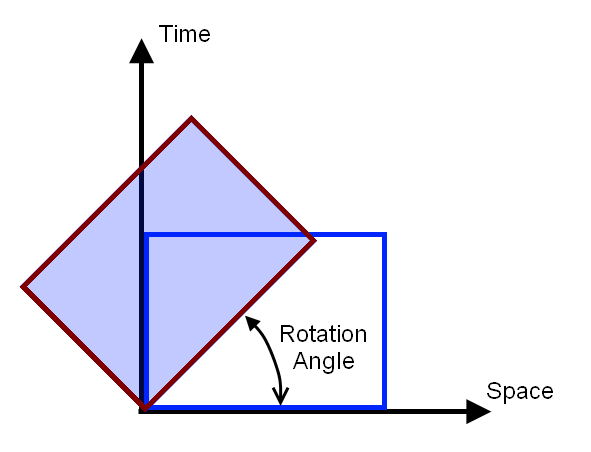}
\caption{Conceptual Diagram of a Domain Rotation in Space-Time}
\label{fig:Rotation}
\end{figure}

To test the effectiveness of other possible domain rotation angles, the domain $\Omega\times[0,T]$ was rotated about the origin counter-clockwise from $\theta = 0^\circ$ to $\theta = 45^\circ$. This is shown conceptually in Figure ~\ref{fig:Rotation} demonstrating how the domain rotates about the origin of the space-time plain. 

After performing this rotation on the domain, a similar initial value problem as the one given in ~\ref{eqn:weakform} was then solved on the new rotated domain. The algebraic formulation rotated domain is given below.

\begin{align}
	\mathit D \widehat \Psi(x',t') = 0 \nonumber \\ 
	\widehat\Psi(\cdot', 0') = \Psi^0 \nonumber \\
	\text{where \;} 
	\begin{bmatrix} x' \\ t'
	\end{bmatrix} = 
	\begin{bmatrix} 
	cos(\theta) & -sin(\theta) \\
	sin(\theta) & cos(\theta) 
		\end{bmatrix} 
		\begin{bmatrix} x \\ t
		\end{bmatrix} 
		\nonumber \\
		\text{and \;} [0 < x < x_{max}] \text{, \;} [0 < t < t_{max}] 	\nonumber
\end{align}

In this test,  $t_{max} = .4$ and $x_{max} = \frac{2}{3}*sec(45^\circ-\theta)$, where $\theta$ is the angle of rotation. The $x$-axis is scaled by a factor of $\sec(45^\circ - \theta)$ so that the wave function remains centered in the domain throughout the rotation. This rotation has the added advantage of following the path of the solution more closely, since high-energy, low mass solutions to the Dirac equation tend to move along the characteristic lines of equation; which is to say, particles that have high energy and low mass move at nearly the speed of light.

This test was also conducted with  non-zero masses, and results compared to a solution calculated using numerical Fourier transformation. 
Due to stability concerns in the massive case, a Dirichlet boundary condition was added to the side $x = 0$ and the the center of the wave function was moved from $x = 0.5$ to $x = 0.8$ to keep the wave function from colliding with the Dirichlet boundary condition on the $x=0$ wall.  The domain shape was slightly altered to $\left [0 \leq x \leq 1.6\right ] \times \left [0 \leq t \leq .4 \right]$.

\begin{figure}[h]
\centering
\subfigure[]{
	\includegraphics[width=0.3\textwidth]{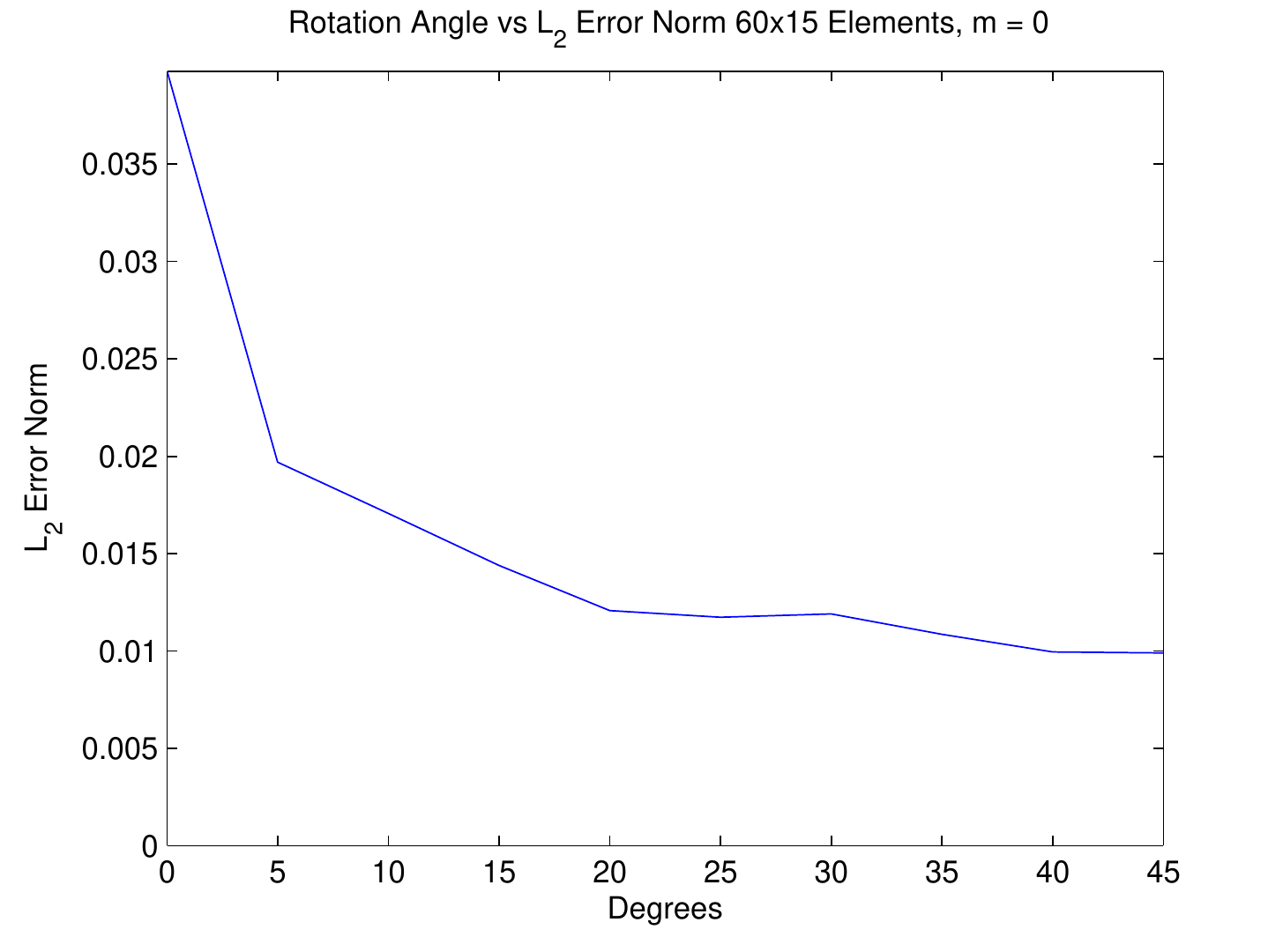} }
\subfigure[]{
	\includegraphics[width=0.3\textwidth]{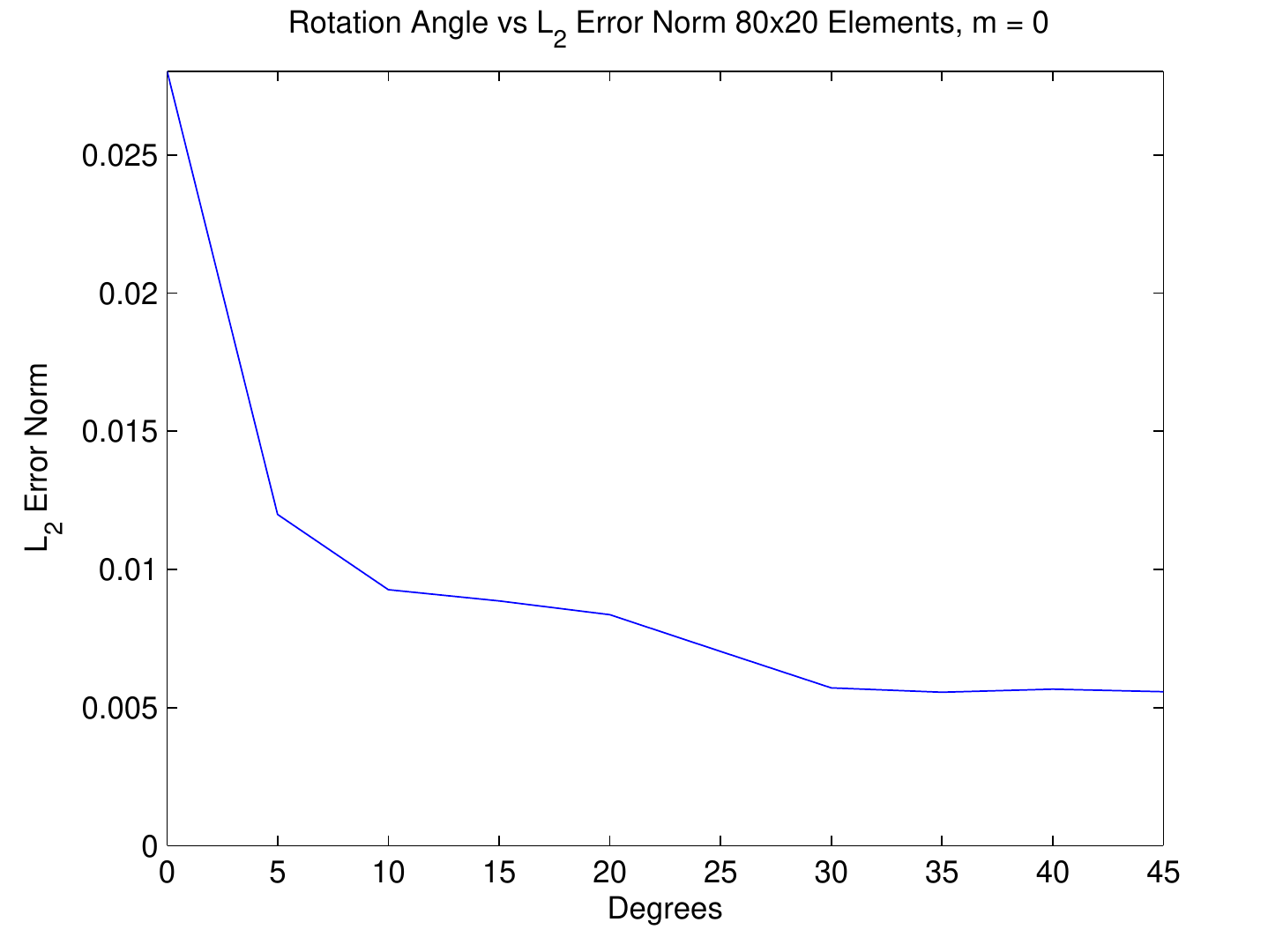} }
\subfigure[]{
	\includegraphics[width=0.3\textwidth]{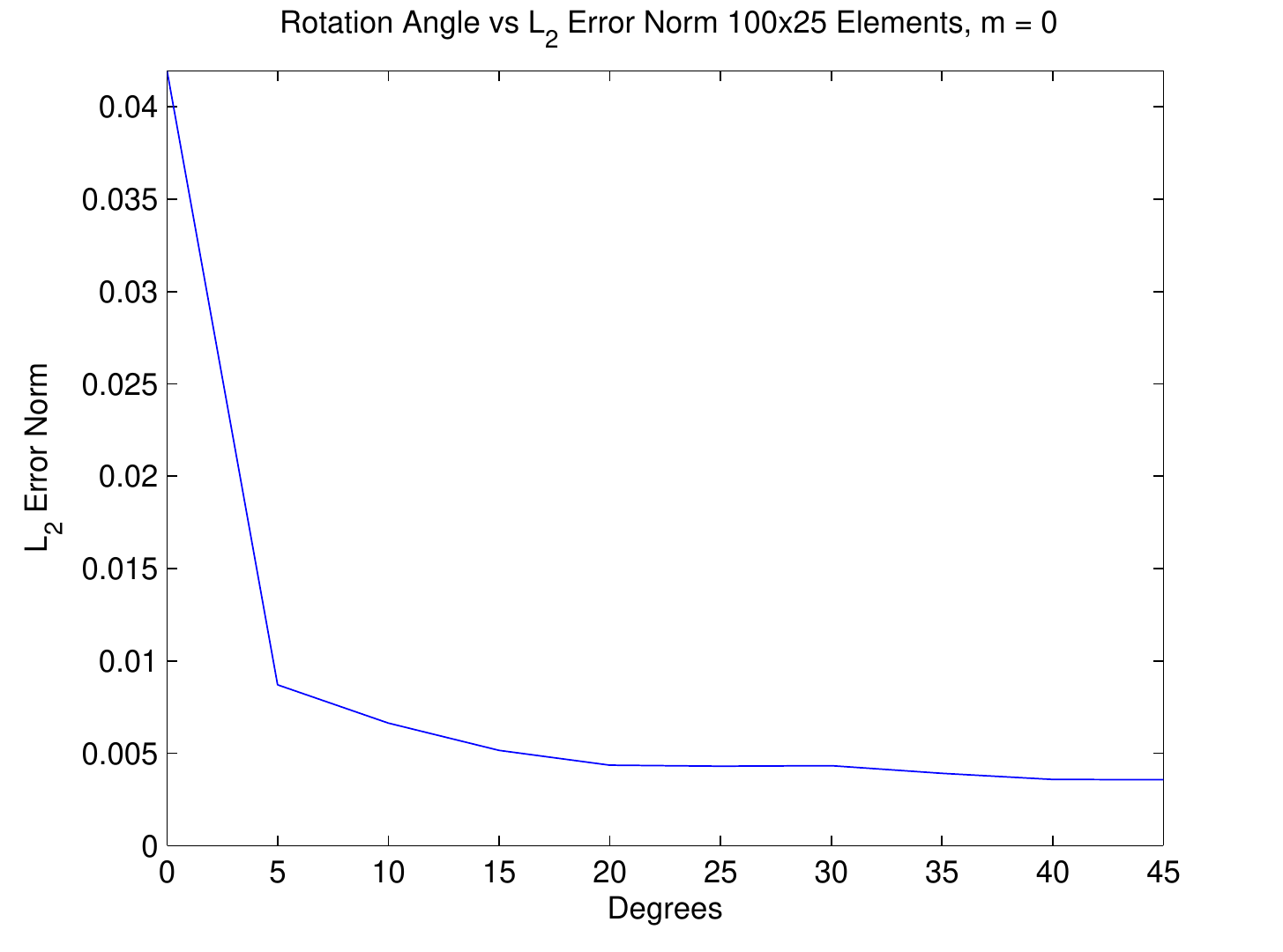} }
\\
\subfigure[]{
	\includegraphics[width=0.3\textwidth]{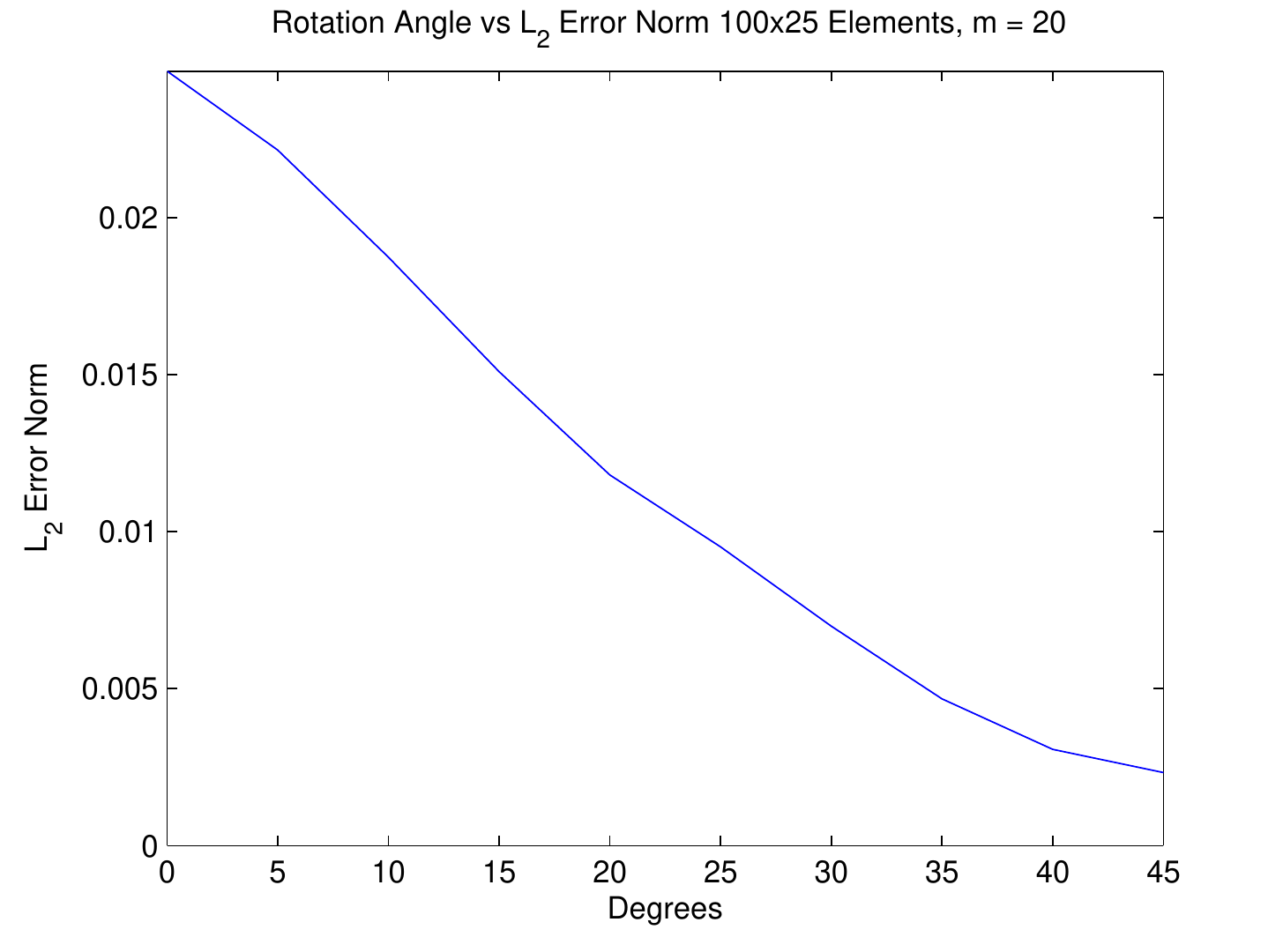} }
\subfigure[]{
	\includegraphics[width=0.3\textwidth]{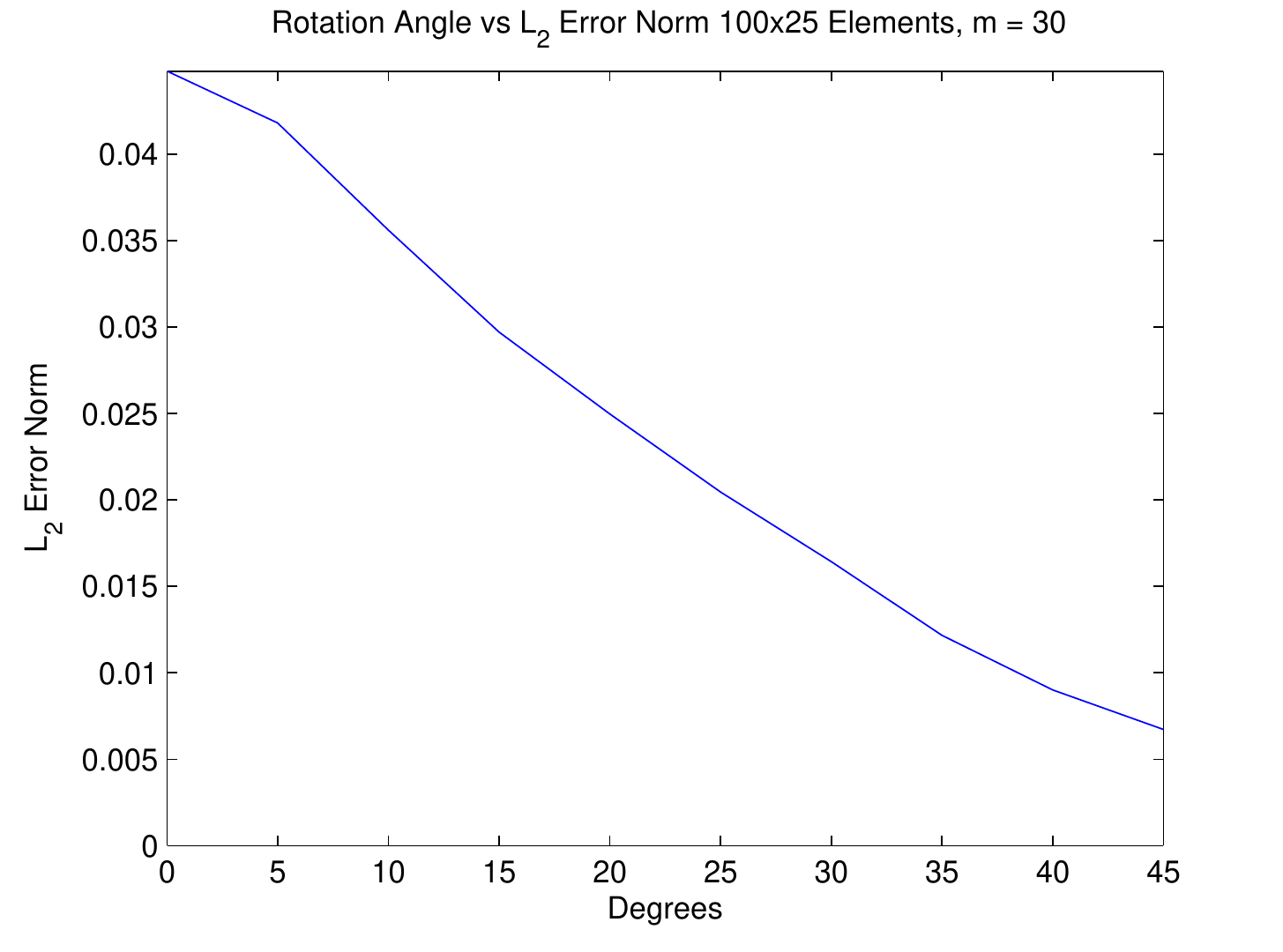} }
\subfigure[]{
	\includegraphics[width=0.3\textwidth]{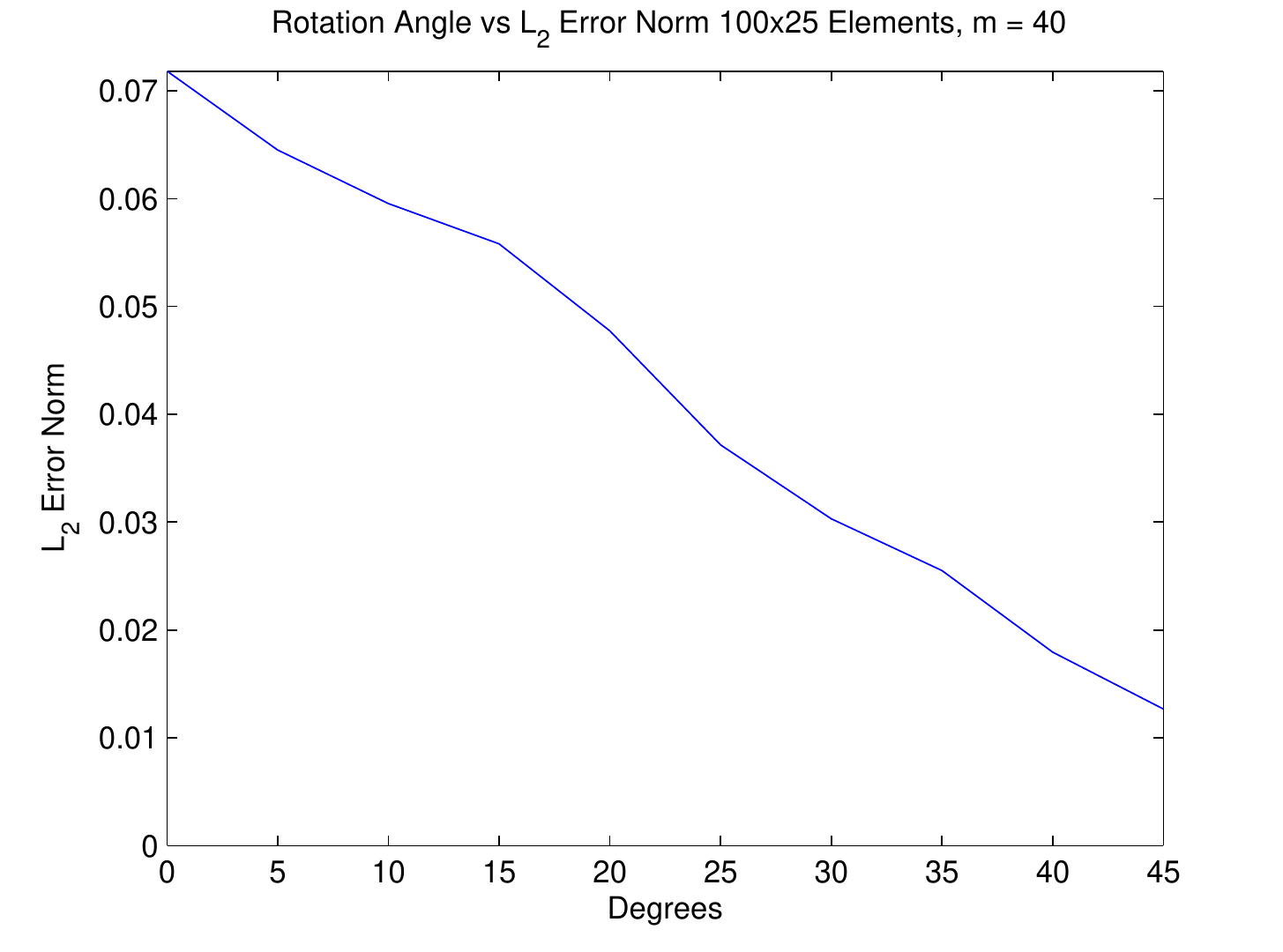} }
\caption{$|\Psi - \Psi_{FEM}|_2$, the $L_2$ Error Norm of the Computed Wave Function vs the Angle of Rotation}
\label{fig:AngleError}
\end{figure}

\begin{figure}[h]
\centering
\subfigure[]{
	\includegraphics[width=0.3\textwidth]{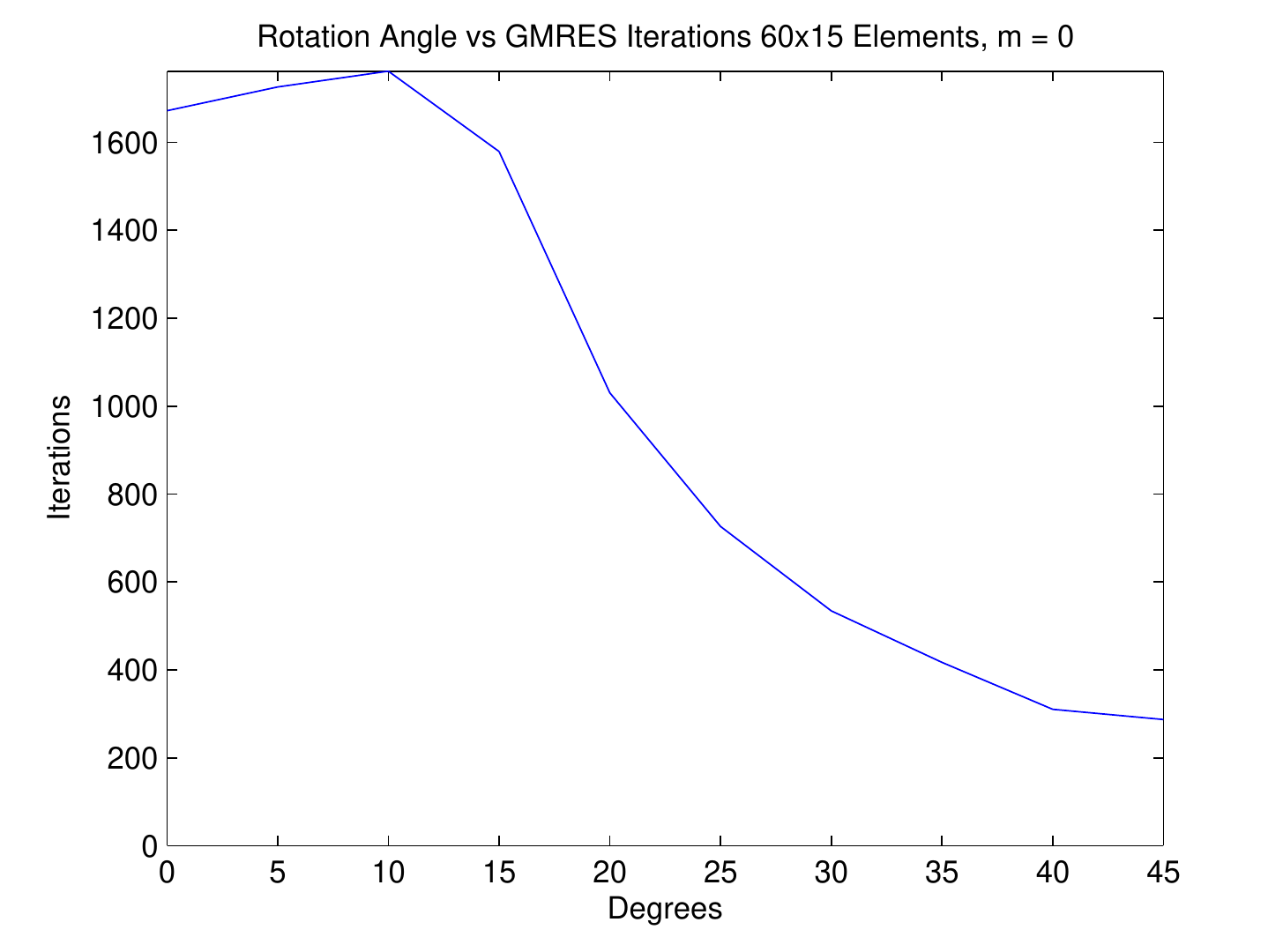}}
\subfigure[]{
	\includegraphics[width=0.3\textwidth]{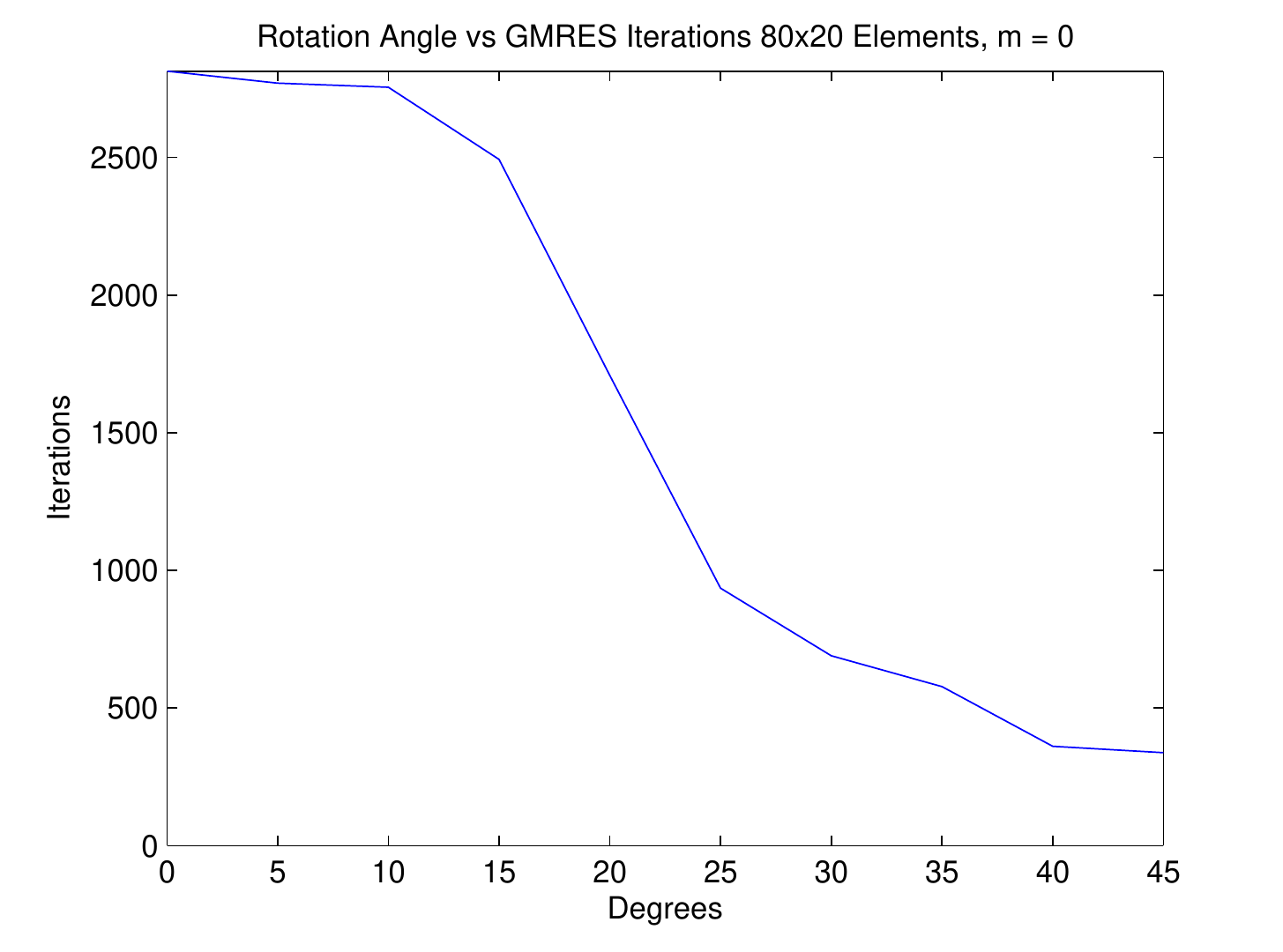}}
\subfigure[]{
	\includegraphics[width=0.3\textwidth]{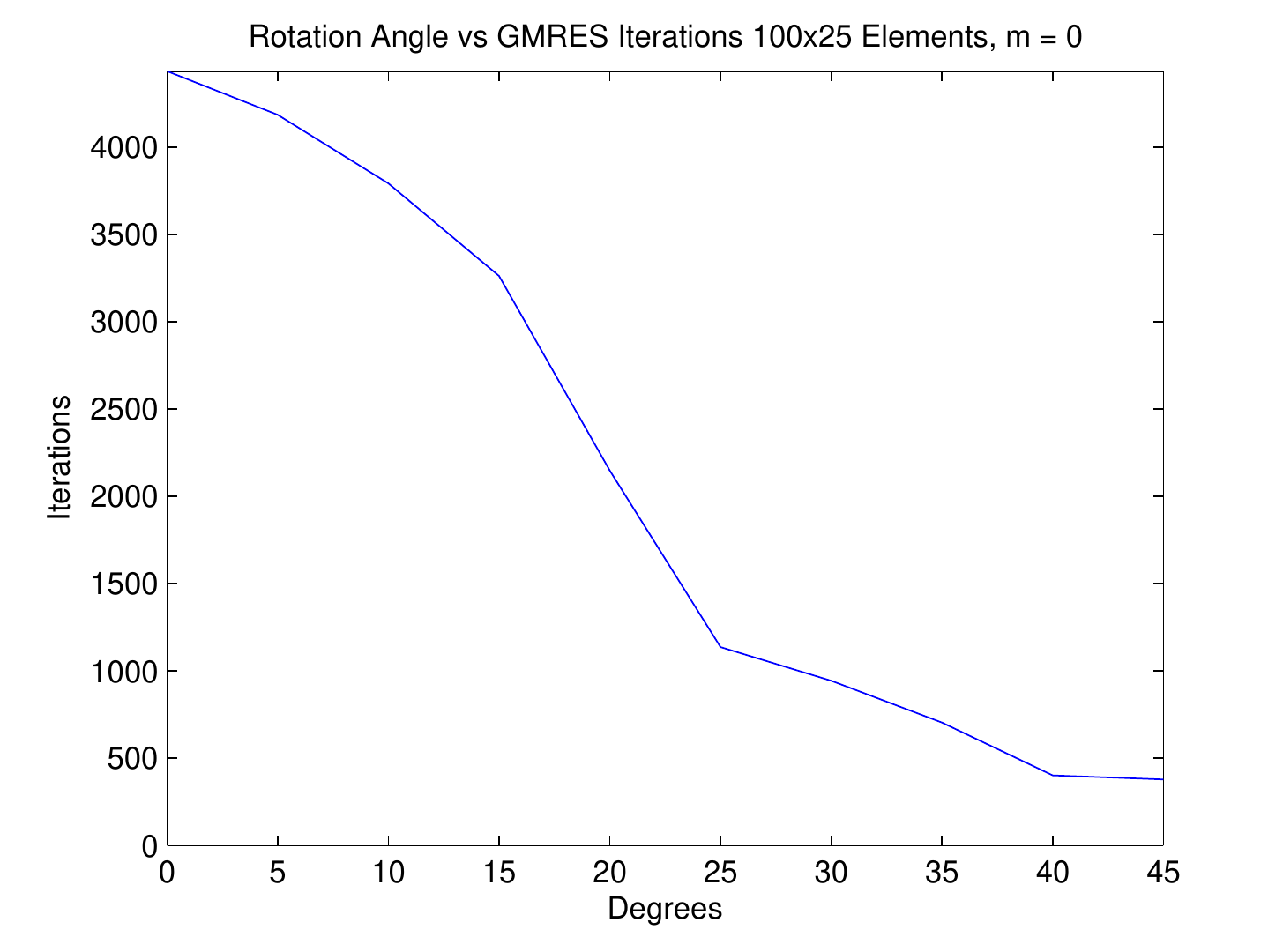}}
\\
\subfigure[]{
	\includegraphics[width=0.3\textwidth]{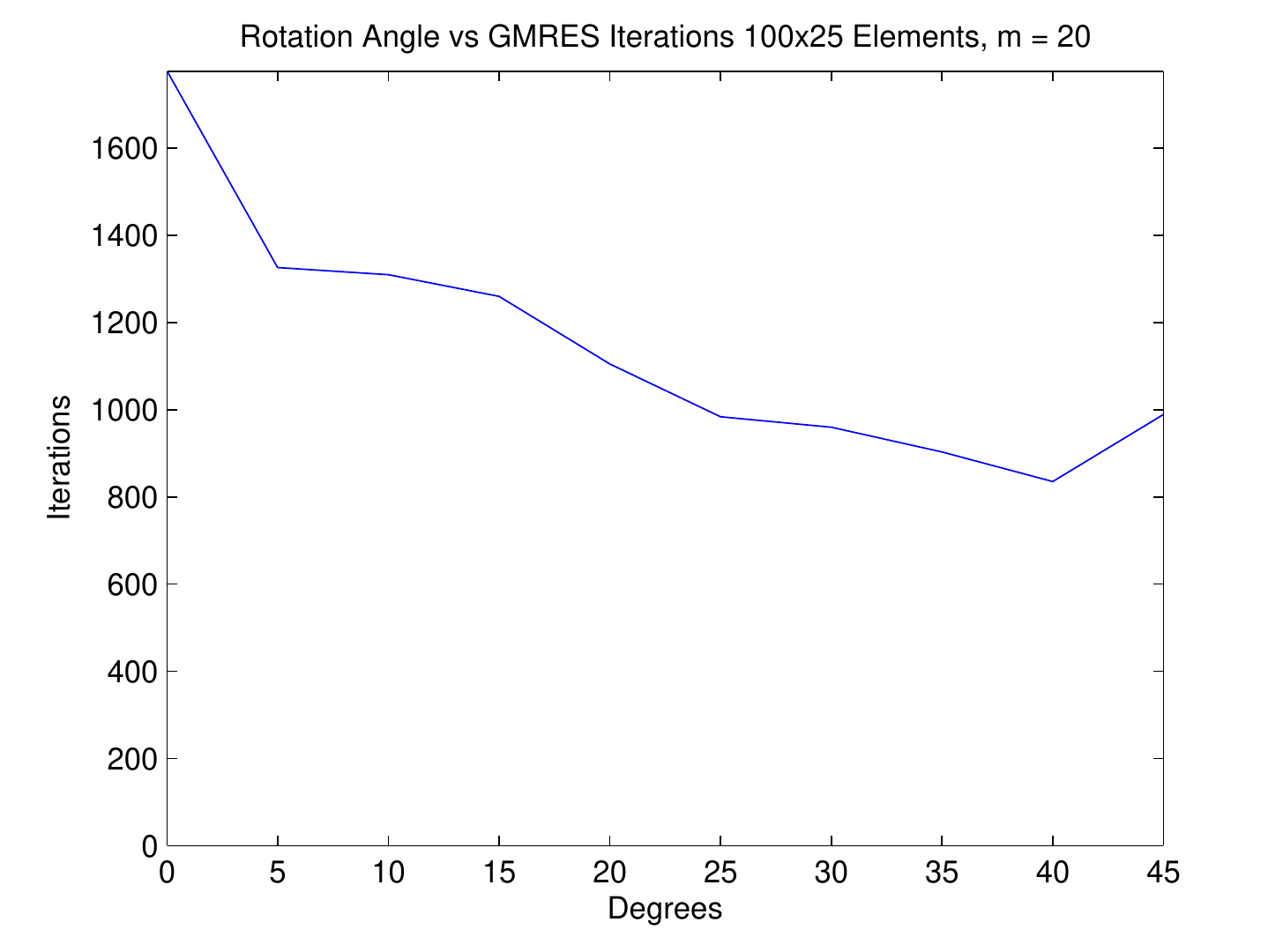}}
\subfigure[]{
	\includegraphics[width=0.3\textwidth]{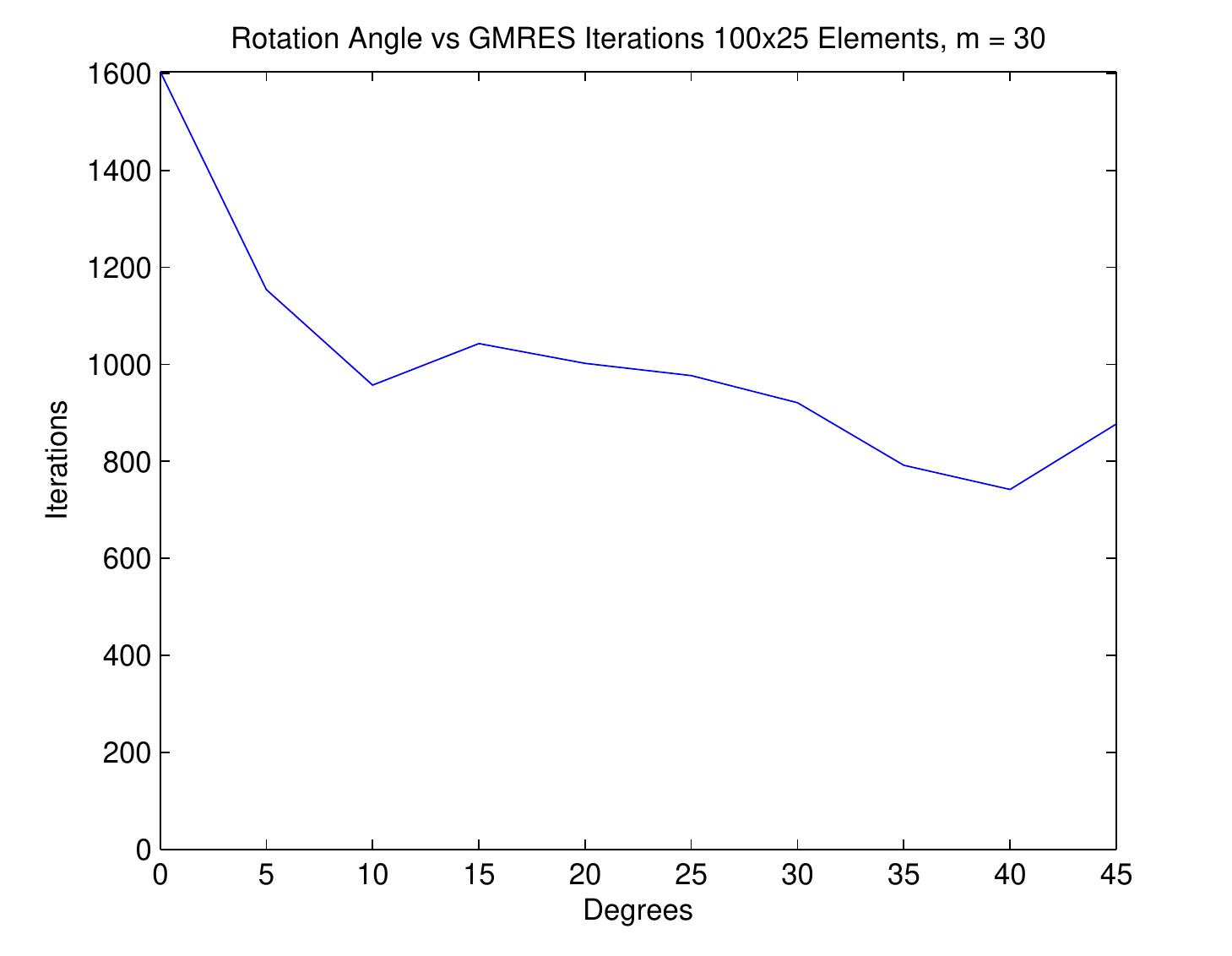}}
\subfigure[]{
	\includegraphics[width=0.3\textwidth]{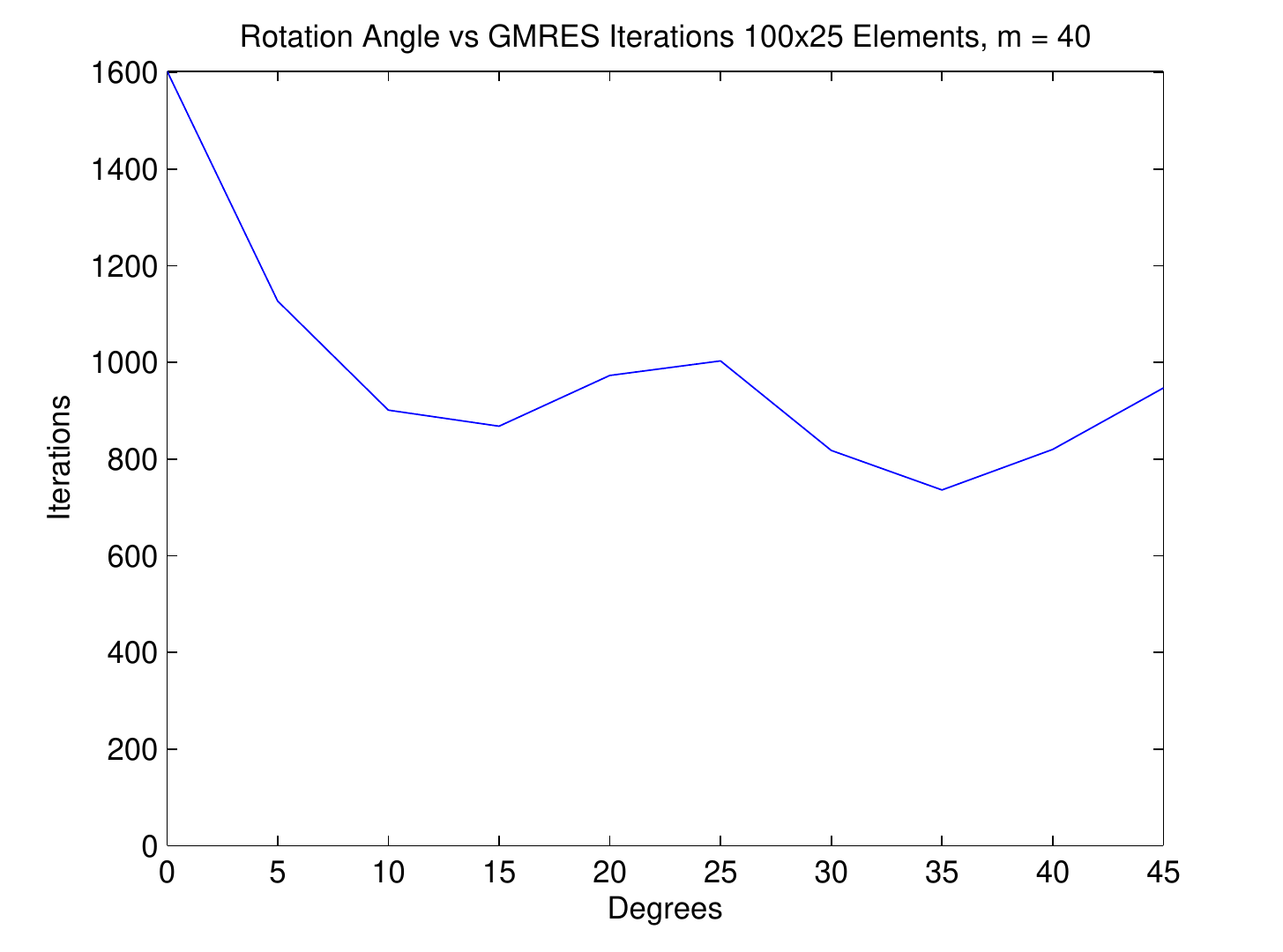}}
\caption{GMRES iterations to reach a residual of $10^{-6}$ vs the Angle of Rotation}
\label{fig:AngleIterations}
\end{figure}

The number of GMRES iterations and the $L_2$ norm of the error were recorded and plotted against the rotation angle used. The results are shown in Figure ~\ref{fig:AngleIterations} and ~\ref{fig:AngleError}. 

%

In Figure ~\ref{fig:AngleError} the $L_2$ norm of the error is shown to decrease steadily as the angle is increased from $0^\circ$ to $45^\circ$ and  reaches its lowest point at $45^\circ$. This relationship is demonstrated in each case tested, regardless of particle mass.

Likewise, in Figure  ~\ref{fig:AngleIterations} the number of GMRES iterations is shown to decrease steadily as the angle is increased from $0^\circ$ to $45^\circ$ and reaches its lowest point at $45^\circ$. This relationship is demonstrated for each domain size tested but only in the massless case.

Interestingly, this relationship is changed somewhat as the mass increases. At a mass of $m = 20$ and $m=30$, the fewest iterations are used at around $40^\circ$. When $m=40$, the fewest iterations are required at $35^\circ$. Finding the source of this small off-angle efficiency improvement requires further investigation.

These results demonstrate that a $45^\circ$ rotation in space-time shows the lowest error for all angles and mass levels tested. Further, they also demonstrate that a $45^\circ$ rotation is either optimal or nearly optimal from a GMRES iterations perspective; however, this relationship is less straightforward than the correlation with error and will require further research to establish the relationship between the massive particle and the optimal rotation angle for algorithmic efficiency.

\section{Conclusion}
\label{sec:Conc}

From the data above we have shown several space-time approaches that may be useful in numerical calculations of the Dirac equation in a rectangular space-time domain. We have also shown that a physically-motivated selected of element and domain shape can substantially improve performance and reduce error for the numerical experiments considered above. Since this improvement was shown for an equation dominated by first-order operators, it may be possible to use a similar approach for other equations with unstable first-order operators as well. 

The results above also show that problems with super-luminal and physically inconsistent propagation may be addressed by choice of discretization and using a fully implicit method. This is then corrected without reference to the problem of Fermion doubling as was suggested by M{\"u}ller et al.  in  \cite{muller1998}.

To expand the usefulness of this numerical approach, further research should be conducted in several areas. One, the sample problem should be developed for $2+1$ and $3+1$ dimensional settings. Two, a more formal error analysis should be conducted to understand the root cause of the simulation behavior above.Three, scalable preconditioners should be investigated for new numerical solutions to the Dirac equation, especially given the size of Dirac-based problems in $3+1$ dimensions (or more). Finally, this model should be tested for suitability and performance in more realistic, inhomogeneous or nonlinear settings.

{\bf Acknowledgements} We thank you to Professor Dongming Mei from the Department of Physics of University of South Dakota for his discussion on  physics interpretation on this results.

\bibliography{P0Bib}
\bibliographystyle{unsrt}

\appendix
\section{Solution of the Initial Value Problem}

This discussion follows closely to the derivation presented in the appendices of \cite{fillion2012}.
In order to compare our results with known solutions of the Dirac equation, we will first consider the case of the massless Dirac equation. 
Since we are interested in the behavior of particles whose mass is very close to zero, this should give us some indication of the fitness of our approach for real-world problems.

Removing the mass term from equation~\ref{eqn:dirac1dOp} and multiplying both sides by the matrix $\BM 1 & 0 \\ 0 & -1 \EM$ gives us the following equation.

\BE
\left( i\hbar \BM 1&0 \\ 0&1 \EM \partial_t  +  i\hbar \BM 0&-1 \\ -1&0 \EM \partial_x  \right) \BM \Psi_l(x,t) \\ \Psi_r(x,t)  \EM = 0 \nonumber 
\EE

We then make the following substitutions 
\BEA
\BM 1&0 \\ 0&1 \EM = I & \BM 0&-1 \\ -1&0 \EM = -\sigma_1 & \BM \Psi_l(x,t) \\ \Psi_r(x,t)  \EM = \Psi(x,\tau) \nonumber
\EEA

and rearrange the equation as follows.

\begin{eqnarray}
-i\hbar I\partial_t \Psi(x,\tau) = - i\hbar c \sigma_1\partial_x  \Psi(x,\tau)  \label{eqn:masslessdirac}
\end{eqnarray}

We may further simplify this equation into a first order ODE by taking the Fourier transform with respect to $\hbar\omega = p$ which expresses the massless Dirac operator in momentum space.

\begin{eqnarray}
-i\hbar I\partial_t \Psi(p,\tau)  =  \sigma_1 p_x  \Psi(p,\tau)  \label{eqn:masslessdiracfourier}
\end{eqnarray}

Integrating directly from $\tau = 0$ to $\tau = t$ gives us the general solution to the massless initial value problem in momentum space.

\begin{eqnarray}
 \Psi(p,t)  =  e^{(\frac i \hbar \sigma_1 p_x t)}  \Psi_0(p)  \nonumber
\end{eqnarray}

Where we may then apply Euler's Identity in order to remove the matrix from the exponential

\begin{eqnarray}
 \Psi(p,t)  =  \left( I cos(\frac {p_x} {\hbar} t) + i \sigma_1 sin(\frac {p_x} {\hbar} t) \right)   \Psi_0(p)  \label{eqn:masslessdiracsolutionPspace}
\end{eqnarray}

Where $I$ is the $2\times2$ identitity matrix. If we take our initial function to be a Gaussian wave  of the form $\Psi_(x,0) = 
\begin{bmatrix}
e^{(i \pi b x - (ax)^2)} \\
e^{(i \pi b x - (ax)^2)} 
\end{bmatrix}$
which may be expressed in momentum space as $\widehat \Psi_0(p)  = 
\begin{bmatrix}
(2a^2)^{-\frac{1}{2}}e^{-\frac{(\omega+ \pi b)^2}{4a^2} } \\
(2a^2)^{-\frac{1}{2}}e^{-\frac{(\omega+ \pi b)^2}{4a^2} }
\end{bmatrix}$.
Inverting the Fourier transform from equation~\ref{eqn:masslessdiracsolutionPspace} with the given initial value results in the general initial value solution:

\begin{equation}
\Psi(x,t)  = 
\begin{bmatrix}
A &  B\\
C  &  D
 \label{eqn:masslessdiracsolution}
\end{bmatrix}
\begin{bmatrix}
\Psi_1(x,0) \\
\Psi_2(x,0)
\end{bmatrix}
\end{equation}

where the values $A$, $B$, $C$, and $D$ are define as

\BEA
A = \frac{1}{2} (e^{ -a^2 (t^2-2xt)-i b \pi t}+e^{-a^2 (t^2+2xt) + i b \pi t}) \nonumber\\
B = \frac{1}{2}( -e^{ -a^2 (t^2-2xt) - i b \pi t}+e^{-a^2 (t^2+2xt) + i b \pi t}) \nonumber\\
C = \frac{1}{2}( -e^{ -a^2 (t^2-2xt) - i b \pi t}+e^{-a^2 (t^2+2xt) + i b \pi t}) \nonumber\\
D = \frac{1}{2} (e^{ -a^2 (t^2-2xt) - i b \pi t}+e^{-a^2 (t^2+2xt) + i b \pi t})  \nonumber
\EEA

Equation~\ref{eqn:masslessdiracsolution} may then be used to calculate the analytic solution to any combination of massless Gaussian wave packets with the packet width given by $a$ and the momentum set by $b$. 

\end{document}